\documentclass[11pt,reqno,a4paper]{amsart}
\usepackage{fancyhdr} 
\usepackage{tikz}
\usetikzlibrary{cd}
\usepackage{amsmath}

\usepackage[active]{srcltx}
\usepackage[pagebackref, colorlinks, linkcolor=red, citecolor=blue, urlcolor=blue, hypertexnames=true]{hyperref}
\usepackage[backrefs]{amsrefs}
\usepackage{enumerate}
\usepackage[mathscr]{eucal}
\allowdisplaybreaks
\usepackage[all,cmtip]{xy}

\usepackage{color}

\allowdisplaybreaks
\usepackage[all]{xy}

\newcommand{\R}{{\mathbb R}}

\newcommand{\Z}{{\mathbb Z}}
\newcommand{\N}{{\mathbb N}}

\def\v{\varepsilon}

\def\sub{\subset}

\theoremstyle{plain}
\numberwithin{equation}{section}

\newtheorem{thm}{Theorem}[section]
\newtheorem{conjecture}[thm]{Conjecture}
\newtheorem{cor}[thm]{Corollary}
\newtheorem{lem}[thm]{Lemma}
\newtheorem{condition}[thm]{Condition}

\newtheorem{definition}[thm]{Definition}
\newtheorem{example}[thm]{Example}
\newtheorem{remark}[thm]{Remark}
\newtheorem{pro}[thm]{Proposition}

\newtheorem{claim}[thm]{Claim}

\theoremstyle{definition}

\subjclass[2010]{28A80; 37B05; 54F45}
\keywords{Amenable group; Minkowski dimension; Continuity; Mean Hausdorff dimension; Metric mean dimension; Product formula}

\begin{document}
	\title{amenable metric mean dimension and amenable mean Hausdorff dimension of product sets and metric varying}

\email{lixq233@mail2.sysu.edu.cn }
\email{luoxf29@mail2.sysu.edu.cn }

\thanks{{$^{*}$}Corresponding author: Xiaofang Luo }

\maketitle
\centerline{\scshape Xianqiang Li, Xiaofang Luo$^*$}

	\maketitle
	\medskip
	\noindent{\bf Abstract} Metric mean dimension and mean Hausdorff dimension depend on metrics. In this paper, we investigate the continuity of the metric mean dimension and mean Hausdorff dimension concerning the metrics for amenable group actions, which extends recent results by Muentes, Becker, Baraviera et al.. Moreover, we give proof of the product formulas for the mean Hausdorff dimension and the metric mean dimension for amenable group actions.

\section{Introduction}
\hspace{4mm}
Mean dimension, introduced by Gromov in 1999 \cite{MG}, is a topological invariant that proves particularly useful for studying the complexity of systems with infinite entropy. Lindenstrauss and Weiss \cite{LW} later introduced the concept of metric mean dimension, which is not independent of the metric and closely related to the ``mean'' Minkowski dimension. They demonstrated that the metric mean dimension serves as an upper bound for the mean dimension across any metric of the given system. Lindenstrauss \cite{LE} showed that, with the marker property, a metric exists where the mean dimension equals the lower metric mean dimension. Similarly, Lindenstrauss and Tsukamoto \cite{LW} proposed that, with the marker property, there is a metric where the mean dimension coincides with the upper metric mean dimension. Additionally, Lindenstrauss and Tsukamoto \cite{LW} introduced a metric-independent concept called the mean Hausdorff dimension, which provides a more precise upper bound for the mean dimension. For more kinds of mean dimensions, one may refer to \cites{CDL, LF,LB,WY,WC1,TM2}. Recently, Muentes, Becker, Baraviera et al. \cite{ABB} have investigated the continuity of both metric mean dimension and mean Hausdorff dimension. We will briefly review their findings.

\hspace{4mm}
Suppose $(X, \tau)$ is a compact topological space and $f :X \to X$ is a continuous map, denote by $X(\tau)$ the set of all metrics that induce the same topology  $\tau$ on $X$, and denote by $\operatorname{mdim}_{\mathrm{H}}(\cdot)$ and $\operatorname{mdim}_{\mathrm{M}}(\cdot)$ the mean Hausdorff dimension and metric mean dimension, respectively. Let $$\mathcal{A}_{d}(X)=\left\{g_{d}: g_{d}(x, y)=g(d(x, y)) \text { for all } x, y \in X\text {, and } g \in \mathcal{A}[0, \rho]\right\} \text {, }$$
where $\rho$  is the diameter of  $X$  and
$$\mathcal{A}[0, \rho]=\left\{g:[0, \rho] \rightarrow[0, \infty): g \text { is continuous, increasing, subadditive and } g^{-1}(0)=\{0\}\right\} .$$
For any continuous map $\zeta\in \mathcal{A}[0, \rho]$, take
$$k_m(\zeta) = \liminf_{\varepsilon \to 0^+} \frac{\log(\zeta(\varepsilon))}{\log(\varepsilon)},\quad k_M(\zeta) = \limsup_{\varepsilon \to 0^+} \frac{\log(\zeta(\varepsilon))}{\log(\varepsilon)}.$$ 

Set
 $${\mathcal{A}}^{+}_d(X) =\{\zeta\circ d \in \mathcal{A}_d(X): \zeta\in {\mathcal{A}}^{+}[0,\rho]   \},$$ 
where
$${\mathcal{A}}^{+}[0,\rho]:=\{\zeta\in {\mathcal{A}}[0,\rho]: k_m(\zeta)=k_M(\zeta)>0\}.$$

\begin{thm}\cite{ABB}
Let  $P=M  \text{ or } H$, $(X, \tau)$ is a compact topological space and $f :X \to X$ is a continuous map, if there is $ d \in X(\tau)$ such that $\operatorname{mdim}_{\mathrm{P}}(X, d, f)>0,$ then
$$\begin{aligned}
\operatorname{mdim}_{\mathrm{P}}(X, f): X(\tau) & \rightarrow \mathbb{R} \cup\{\infty\} \\
d & \mapsto \operatorname{mdim}_{\mathrm{P}}(X, d, f)
\end{aligned}$$
is not continuous anywhere.
\end{thm}
\begin{thm}\cite{ABB}
Let  $X$  be a compact space such that $d: X \times X \rightarrow[0, \rho]$ is a surjective metric map. If $\operatorname{mdim}_{\emph{M}}(X, f, d)<\infty,$ then the maps
$$\begin{aligned}
\overline{\operatorname{mdim}}_{\mathrm{M}}(X, f):\left(\mathcal{A}_{d}^{+}(X), \mathcal{W}\right) & \rightarrow \mathbb{R} \\
g_{d} & \mapsto \overline{\operatorname{mdim}}_{\mathrm{M}}\left(X, g_{d}, f\right)
\end{aligned}$$
and
$$\begin{aligned}
\underline{\operatorname{mdim}}_{\mathrm{M}}(X, f):\left(\mathcal{A}_{d}^{+}(X), \mathcal{W}\right) & \rightarrow \mathbb{R} \\
g_{d} & \mapsto \underline{\operatorname{mdim}}_{\mathrm{M}}\left(X, g_{d}, f\right)
\end{aligned}$$
are continuous.
\end{thm}

\hspace{4mm}
Many classical theories of mean dimension have been extended to a broader range of group actions, including $\mathbb{Z}^k$-actions \cites{DJM, GLT}, amenable group actions \cites{CDZ,LW,KF,LL}, and sofic group actions \cites{HL,YG,HB,LLB}. This naturally raises the question of whether the continuity of metric mean dimension and mean Hausdorff dimension can be extended to amenable group actions. In this paper, we explore these dimensions within the framework of amenable group actions, focusing on their continuity concerning metrics. Our main results build upon the methods of Muentes, Becker, Baraviera et al., but some of our conclusions are a bit different, since the metric mean dimension and mean Hausdorff dimension may not always exist for any given dynamical system, we consider the upper and lower metric mean dimension and upper and lower mean Hausdorff dimension for amenable group actions, respectively.

\hspace{4mm}
In this paper, we also concern the product formulas for the mean Hausdorff dimension and the metric mean dimension for amenable group actions. For mean dimension $\operatorname{mdim}(\cdot)$, the product formula 
$$\operatorname{mdim}(X\times Y)\leq \operatorname{mdim}(X)+ \operatorname{mdim}(Y)$$
is familiar to us \cite{LW}. In 2019, Tsukamoto \cite{TM} provided examples demonstrating that this inequality can be strictly. Recently, Jin and Qiao \cite{JQ} examined the inequality in the case where $X=Y$, and derived an interesting formula for the mean dimension of product spaces. Liu, Selmi, and Li \cite{LSL} investigated product formulas for the mean Hausdorff dimension and the metric mean dimension, presenting one of their key results in \cite{LSL} Theorem 3.21. However, a crucial part of their proof, specifically Lemma 3.19, requires improvements, which affects the proofs of Lemma 3.20 and Theorem 3.21 in \cite{LSL}. In this paper, we provide our proof (see Lemma \ref{2.11} to Theorem \ref{a2.2}) within the framework of amenable groups. 

\hspace{4mm}
The organization of this paper is as follows. Section \ref{sec2} reviews key concepts related to amenable group actions, including metric mean dimension, mean Hausdorff dimension, and Katok entropy, and explores their interrelationships. Section \ref{sec3} establishes product formulas for the metric mean dimension and mean Hausdorff dimension within the framework of amenable group actions, providing an illustrative example and deriving formulas for the Minkowski dimension and metric mean dimension. Section \ref{sec4} examines the continuity of metric mean dimension and mean Hausdorff dimension concerning metrics for amenable group actions. Finally, section \ref{sec5} addresses the continuity of metric mean dimension in specific metric spaces for amenable group actions and presents three examples with detailed expressions of metric mean dimensions for particular metrics.

\section{Preliminaries}\label{sec2}
\hspace{4mm}		 
Let $G$ be a group. One says that a sequence $\{F_{n}\}_{n \geq 1 }$ of non-empty finite subsets of $G$ is a \textbf{F$\phi$lner sequence} for $G$ if one has
  $$ \lim_{n \rightarrow \infty}\frac{|F_{n} \setminus gF_{n}|}{|F_{n}|}=0.
  $$

  \hspace{4mm}
We say that a countable group is \textbf{amenable} if it admits a F$\phi$lner sequence. A F$\phi$lner sequence $\{F_{n}\}$ in $G$ is said to be \textbf{tempered} if there exists a constant $C>0$ which is independent of $n$ such that 
$$|\bigcup_{k <n}F_{k}^{-1}F_{n}|\leq C|F_{n}|, \text{for every } n\in \N.
$$

\begin{definition}
     Let $G$ be a countable discrete amenable group. By a pair $(X,G)$ we mean a \textbf{$G$-system}, where $X$ is a compact metric space and $\Gamma:G \times X \rightarrow X, \text{ given by } (g,x) \rightarrow gx$,  is a continuous mapping satisfying:\\
(1)$\Gamma(1_{G},x)=x$ for every $x \in X$;
\\(2)$\Gamma(g_{1},\Gamma(g_{2},x))=\Gamma(g_{1}g_{2},x)$ for every $g_{1},g_{2} \in G$ and $x\in X$.
\end{definition}

\hspace{4mm}
In this paper, we always assume $G$ is a countable discrete amenable group and $X$ is a compact metric space. Let $\text{Fin}(G)$ be the family of finite nonempty subsets of $G$.

\hspace{4mm}
Let $(X,G)$ be a $G$-system with a metric $d$. For $F \in \text{Fin}(G)$, define a metric $d_{F}$ on $X$ by
$$d_{F}(x,y)=\max_{g \in F}d(gx,gy), \text{~for~every~} x, y \in X.
$$

\begin{definition}
   Let $E$ be a compact subset of $X$, for any $F \in \emph{Fin}(G)$ and $\varepsilon>0$, let $K \sub  E $, if for any $x \in E$ there exists $y \in K$ such that $d_{F}(x,y) \leq \varepsilon$, then we call $K$ an $(F,\varepsilon)$\textbf{-spanning set} of $E$. $K$ is called an $(F,\varepsilon)$\textbf{-separated set} if we have $d_{F}(x,y) >\varepsilon$, for any distinct $x,y \in K$.
\end{definition} 

\hspace{4mm}
Denote by $s_{F}(d,\varepsilon,E)$ the maximal cardinality of any $(F,\varepsilon)$-separated subset of $E$, by $r_{F}(d,\varepsilon,E)$ the smallest cardinality of any $(F,\varepsilon)$-spanning subset of $E$, by $\text{cov}_{F}(d,\varepsilon,E)$ the smallest cardinality of any open cover $\alpha$ of $E$ that satisfies $\text{mesh}(\alpha,d_{F}) <\varepsilon$, where $\text{mesh}(\alpha,d_{F}):=\max_{A \in \alpha}\text{diam}_{d_{F}}(A)$. 

\hspace{4mm}
Let $\{F_{n}\}$ be any F$\phi$lner sequence in $G$, denote
\begin{itemize} \item $s(d,\varepsilon,\{F_{n}\},E)=\underset{n\to\infty}\limsup\frac{1}{|F_{n}|}\log  s_{F_{n}}(d,\varepsilon,E)
$;
\item $r(d,\varepsilon,\{F_{n}\},E)=\underset{n\to\infty}\limsup\frac{1}{|F_{n}|}\log
 r_{F_{n}}(d,\varepsilon,E)$;
\item  $\text{cov}(d,\varepsilon,\{F_{n}\},E)=\underset{n\to\infty}\limsup \frac{1}{|F_{n}|}\log \text{cov}_{F_{n}}(d,\varepsilon,E)$.
\end{itemize}   

\begin{remark}\cite{LZ}\label{01010}
Let $(X,G)$ be a $G$-system with a metric $d$, $\{F_{n}\}$ a F$\phi$lner sequence in $G$, then for any $\varepsilon>0$ and any compact set $E \sub X$, we have 
\begin{itemize} 
    \item[(i)] $r_{F_{n}}(d,\varepsilon,E) \leq s_{F_{n}}(d,\varepsilon,E) \leq 
\emph{cov}_{F_{n}}(d,\varepsilon,E)$.
\item[(ii)]$\underset{n\to\infty}\lim \frac{1}{|F_{n}|}\log \emph{cov}_{F_{n}}(d,\varepsilon, E)$ always exists and does not depend on the choice of the F$\phi$lner sequence $\{F_{n}\}$.
\end{itemize}

\end{remark}

\begin{definition}\cite{LZ}
  Let $(X,G)$ be a $G$-system, $\{F_{n}\}$ a F$\phi$lner sequence, for any compact set $E \sub X$, the \textbf{upper and lower metric mean dimension} of $E$ are defined by
 $$
\overline{\emph{mdim}} _ {\emph{M}} (E,G,d)=\limsup_{\varepsilon \rightarrow 0 }\frac{s(d,\varepsilon,\{F_{n}\},E)}{|\log \varepsilon|}=\limsup_{\varepsilon \rightarrow 0 }\frac{r(d,\varepsilon,\{F_{n}\},E)}{|\log \varepsilon|}=\limsup_{\varepsilon \rightarrow 0 }\frac{\emph{cov}(d,\varepsilon,\{F_{n}\},E)}{|\log \varepsilon|},
 $$
 $$
\underline{\emph{mdim}} _ {\emph{M}} (E,G,d)=\liminf_{\varepsilon \rightarrow 0 }\frac{s(d,\varepsilon,\{F_{n}\},E)}{|\log \varepsilon|}=\liminf_{\varepsilon \rightarrow 0 }\frac{r(d,\varepsilon,\{F_{n}\},E)}{|\log \varepsilon|}=\liminf_{\varepsilon \rightarrow 0 }\frac{\emph{cov}(d,\varepsilon,\{F_{n}\},E)}{|\log \varepsilon|}.
 $$

These values do not depend on the choice of the F$\phi$lner sequence $\{F_{n}\}$. When the above two values coincide,
	it is called the \textbf{metric mean  dimension} of $E$ and $\mathrm{denoted~by~dim}_{\mathrm{M}}(E,G,d).$ 
\end{definition}
 \hspace{4mm}
Now we recall the definition of the Hausdorff dimension. 
Let $E$ be a compact subset of $X$, for $s\geq 0$ and $\varepsilon> 0$, we define $\text{H} _\varepsilon^s( E, d)$ by
	$$\text{H}_{\varepsilon}^{s}(E,d)=\inf\left\{\sum_{i=1}^{\infty}\left(\operatorname{diam}E_{i}\right)^{s}\Bigg|E=\bigcup_{i=1}^{\infty}E_{i}\text{ with diam }E_{i}<\varepsilon(\forall i\geq1)\Bigg\}.\right. $$
	
By  convention we consider $0^{0}=1$ and $\text{diam}(\emptyset)^{s}=0$. Let $\varphi >0$, take $$ \text{dim}_{\text{H}}(E,d,\varepsilon,\varphi)=\sup\{s \geq 0:  \text{H}_{\varepsilon}^{s} (E,d) \geq \varphi\}.  $$
And set $$ \text{dim}_{\text{H}}(E,d,\varepsilon ):= \text{dim}_{\text{H}}(E,d,\varepsilon,1).$$
Then the \textbf{Hausdorff dimension} 
is given by   $$ \text{dim}_{\text{H}}(E,d):= \lim_{\varepsilon\rightarrow 0}\text{dim}_{\text{H}}(E,d,\varepsilon ).  $$
 
The Hausdorff dimension has an equivalent definition. If we set   
$$\text{H}^s(E,d)=\lim_{\varepsilon \rightarrow 0}\text{H}^s_\varepsilon(E,d),$$ then \textbf{Hausdorff dimension$^*$} \cite{FK}, denoted by $\dim_{\text{H}}^\ast(E,d)$, is given by
$$\dim_{\text{H}}^\ast(E,d)=\inf\{s\geq 0: \text{H}^s(E,d)=0\}=\sup\{s\geq 0: \text{H}^s(E,d)=\infty\}.$$
	
\begin{remark}\label{8999}
    Fix any $\varphi >0$, 
Muentes \cite{ABB} proved that $$\dim_{\emph{H}}(E,d)=\dim_{\emph{H}}^\ast(E,d)$$ and furthermore, $$ \dim_{\emph{H}}^{\varphi}(E,d):=\lim_{\varepsilon\rightarrow 0}\emph{dim}_{\emph{H}}(E,d,\varepsilon,\varphi)=\dim_{\emph{H}}(E,d).$$
\end{remark}

\begin{definition}
  Let $(X,G)$ be a $G$-system, let $\{F_{n}\}$ be a F$\phi$lner sequence in $G$, for any compact set $E \sub X$, the \textbf{upper and lower mean Hausdorff dimensions} of $E$ with respect to $\{F_{n}\}$ are defined by
	$$\begin{aligned}
		\overline{\operatorname{mdim}}_{\mathrm{H}}(E,\{F_{n}\},d)&=\lim_{\varepsilon\to0}\left(\operatorname*{limsup}_{n\to\infty}\frac{\operatorname{dim}_{\mathrm{H}}\left(E,d_{F_{n}},\varepsilon\right)}{|F_{n}|}\right),\\
		\underline{\operatorname{mdim}}_{\mathrm{H}}(E,\{F_{n}\},d)&=\lim_{\varepsilon\to0}\left(\operatorname*{liminf}_{n\to\infty}\frac{\operatorname{dim}_{\mathrm{H}}\left(E,d_{F_{n}},\varepsilon\right)}{|F_{n}|}\right).
	\end{aligned}$$
	These values depend on the choice of the F$\phi$lner sequence $\{F_{n}\}$. 
 When the above two values coincide, the common value is called the \textbf{mean Hausdorff dimension} of $E$ concerning $\{F_{n}\}$ and denoted by $\mathrm{mdim}_{\mathrm{H}}(E,\{F_{n}\},d).$ 
\end{definition}

\hspace{4mm}
Next, we will give an equivalent definition for the mean Hausdorff dimension, we will use the following lemma.

\begin{lem}\label{1.33}\cite{ABB}
    Suppose that $(X,d)$ is a compact metric space. Let $E$ be a compact subset of $X$,
 for $s\geq 0$ and $\varepsilon >0$, set $$  \emph{B}_{\varepsilon}^{s} (E,d)=\inf\left\{ \Sigma_{n=1}^{m}(\emph{diam} (B_{n}))^{s}\begin{array}{c|c}&\{ B_{n}\}_{n=1}^{m} \emph{  is a   cover of }E \emph{ by open}\\
 &\emph{balls with }\emph{diam} (B_n) \leq  \varepsilon\end{array}\right\}.$$   Set  
 $$   {\emph{dim}}^{\star}_{\emph{H}}(E,d,\varepsilon)=\sup\{s \geq 0:  \emph{B}_{\varepsilon}^{s} (E,d) \geq 1\},  $$
 we have that
 $$\emph{dim}_{\emph{H}}(E,d)=\lim_{\varepsilon \rightarrow 0}\emph{dim}^{\star}_{\emph{H}}(E,d,\varepsilon).
 $$
\end{lem}

 \hspace{4mm}
With the proof of Lemma \ref{1.33}, it is not difficult to have the following definition.
\begin{definition}
  Let $(X,G)$ be a $G$-system, let $\{F_{n}\}$ be a F$\phi$lner sequence in $G$, for any compact set $E \sub X$, the \textbf{upper and lower mean Hausdorff dimensions} of $E$ are defined by
	$$\begin{aligned}
		\overline{\operatorname{mdim}}_{\mathrm{H}}(E,\{F_{n}\},d)&=\lim_{\varepsilon\to0}\left(\operatorname*{limsup}_{n\to\infty}\frac{\operatorname{dim}_{\mathrm{H}}^{\star}\left(E,d_{F_{n}},\varepsilon\right)}{|F_{n}|}\right),\\
		\underline{\operatorname{mdim}}_{\mathrm{H}}(E,\{F_{n}\},d)&=\lim_{\varepsilon\to0}\left(\operatorname*{liminf}_{n\to\infty}\frac{\operatorname{dim}_{\mathrm{H}}^{\star}\left(E,d_{F_{n}},\varepsilon\right)}{|F_{n}|}\right).
	\end{aligned}$$
	These values depend on the choice of the F$\phi$lner sequence $\{F_{n}\}$. 
 When the above two values coincide, the common value is called the \textbf{mean Hausdorff dimension} of $E$ concerning $\{F_{n}\}$ and denoted by $\mathrm{mdim}_{\mathrm{H}}(E,\{F_{n}\},d).$ 
\end{definition}
	
\begin{remark}\label{x000}\cite{LL}
    Let $G$ be a countable discrete amenable group. Let $(X, G)$ be a $G$-system with a metric $d$. Let $\{F_{n}\}$ be a F$\phi$lner sequence in $G$, then for any compact set $E \sub X$, 
    $$ \underline{\mathrm{mdim}}_{\mathrm{H}}(E, \{F_{n}\}, d)\le\overline{\mathrm{mdim}}_{\mathrm{H}}(E,\{F_{n}\},d)      
		\le\underline{\mathrm{mdim}}_{\mathrm{M}}(E,G,d)\le\overline{\mathrm{mdim}}_{\mathrm{M}}(E,G,d).
    $$
\end{remark}

\hspace{4mm}
Now we review Katok entropy for $G$-systems. Let $M(X, G)$ be the collection of $G$-invariant probability measures of $X$ and $E(X, G)$ be the set of ergodic measures. Given $F \in \text{Fin}(G)$, let $0 <\delta <1, \varepsilon>0$ and $\mu \in M(X,G)$, a set $D \sub X$ is said to be an $(F,\varepsilon,\delta)$-spanning set if the union $\bigcup_{x \in D}B_{d_{F}}(x,\varepsilon)$ has $\mu$-measure more than $1-\delta$. Let $r_{F}(d,\mu,\varepsilon,\delta,X)$ denote the minimum cardinality of $(F,\varepsilon,\delta)$-spanning sets. Let $\{F_{n}\}$ be a F$\phi$lner sequence in $G$, define the \textbf{Katok $\varepsilon$-entropy} with respect to $\{F_{n}\}$ by
$$h_{\mu}^{K}(\varepsilon,\delta,\{F_{n}\})=\limsup_{n \rightarrow \infty}\frac{1}{|F_{n}|}\log r_{F_{n}}(d,\mu,\varepsilon,\delta,X).$$

Huang and Liu \cite{HLZ} proved that if $\{F_{n}\}$ is a tempered F$\phi$lner sequence of $G$ with $\frac{|F_{n}|}{\log n} \rightarrow +\infty$, and $\mu$ is an ergodic and $G$-invariant Borel probability measure, then 
$$\lim_{\varepsilon \rightarrow 0}h_{\mu}^{K}(\varepsilon,\delta,\{F_{n}\})=h_{\mu}(X,G)$$
for every $\delta \in (0,1)$, where $h_{\mu}(X,G)$ is the measure theoretic entropy of $\mu$ (for a precise definition, see \cite{HLZ} \cite{LZ} ). 

\hspace{4mm}
For any finite measurable partition $\mathcal{P}$ and $r>0$, let $U_{r}(A):=\{x \in A: \exists y \in A^{c}, \text{~with~} d(x,y) <r\}$ and $U_{r}(\mathcal{P}):=\bigcup_{A 
\in \mathcal{P}}U_{r}(A)$. Since $\bigcap_{r>0}U_{r}(\mathcal{P})= \partial\mathcal{P}$, then we have $\lim_{r \rightarrow0}\mu(U_{r}(\mathcal{P}))=\mu(\partial\mathcal{P})$ for any $\mu \in M(X,G)$, where $\partial\mathcal{P}:=\bigcup_{A \in \mathcal{P}}\partial A$ and $\partial A$ is the boundary of $A$.

\hspace{4mm}
If a finite measurable partition $\mathcal{P}$ satisfies $\mu(\partial \mathcal{P})=0$ for some $\mu \in E(X,G)$, then for any $\gamma>0$, we can find $0 <r< \gamma$ such that $\mu(U_{r}(\mathcal{P}))<\gamma$. Let $r_{\mu,\gamma}:=\sup\{r \in \R^{+}: \exists \text{~finite measurable partition~} \mathcal{P} \text{~with~} \mu(\partial\mathcal{P})=0, \text{diam}(\mathcal{P})<\gamma \text{~and~} \mu(U_{r}(\mathcal{P}))<\gamma\}$ and $r_{\gamma}:=\inf_{\mu \in E(X,G)}r_{\mu,\gamma}.$

\begin{condition}\label{1000}
    For any $\gamma>0, r_{\gamma}>0$ and $\lim_{\gamma \rightarrow0}\frac{\emph{log}r_{\gamma}}{\emph{log}\gamma}=1.$
\end{condition}

\hspace{4mm}
Possible examples of dynamical systems satisfying the above condition are the one dimensional uniquely ergodic systems whose ergodic measure is the Lebesgue measure.

\hspace{4mm}
The following conclusion comes from \cite{LZ}, it shows the relationship between Katok entropy and metric mean dimension.

\begin{thm}\label{1001}
    Let $(X,G)$ be a $G$-system with a metric $d$ satisfying Condition \ref{1000}. For any tempered F$\phi$lner sequence $\{F_{n}\}$ in $G$ with 
    $$\lim_{n \rightarrow \infty}\frac{|F_{n}|}{\emph{log}n}=\infty,$$
    we have
    $$\overline{\emph{mdim}} _ {\emph{M}} (X,G,d)=\limsup_{\varepsilon \rightarrow0}\frac{\sup_{\mu \in M(X,G)}h_{\mu}^{K}(\varepsilon,\delta,\{F_{n}\})}{|\emph{log}\varepsilon|},$$
    $$\underline{\emph{mdim}} _ {\emph{M}} (X,G,d)=\liminf_{\varepsilon \rightarrow0}\frac{\sup_{\mu \in M(X,G)}h_{\mu}^{K}(\varepsilon,\delta,\{F_{n}\})}{|\emph{log}\varepsilon|}.$$
\end{thm}

\section{Some fundamental properties}\label{sec3}
\hspace{4mm}
In this section, We study the product formulas for metric mean dimension and mean Hausdorff dimension for amenable group actions, respectively, and we obtain formulas for metric mean dimension and Minkowski dimension.

\hspace{4mm}
Let $G$ be a countable discrete amenable group that acts continuously on compact metric spaces $X$ and $Y$, the \textbf{product action} of $G$ on the product space $X \times Y$ is defined as follows:
$$g(x,y)=(gx,gy), \text{~for all~} g \in G, (x,y) \in X\times Y.
$$

\hspace{4mm}
Let $(X,G)$ and $(Y,G)$ be two $G$-systems, where $(X,d)$ and $(Y,d')$ are compact metric spaces with metrics $d$ and $ d^{\prime}$, respectively. We will endow the product space $X \times Y$ with the metric 
\begin{equation*}
(d\times d^{\prime})((x_{1},y_{1}),(x_{2},y_{2}))=\max\{d(x_{1},x_{2}) , d^{\prime}(y_{1},y_{2})\}, \text{ for } x_1,x_2 \in X \text{ and } y_1,y_2 \in Y.
\end{equation*} 

\hspace{4mm}
First of all, we consider the product formula for metric mean dimension for amenable group actions. We start with the following lemma which is important for our results.

\begin{lem}\label{2.11}
    Let $(X,G)$ and $(Y,G)$ be two $G$-systems, where $(X,d)$ and $(Y,d')$ are compact metric spaces. If $M \sub X$ and $L \sub Y$ are compact sets, then for any $\varepsilon>0$ and $F \in \emph{Fin}(G)$, we have
    $$r_{F}(d\times d^{\prime},\varepsilon,M \times L) \leq r_{F}(d,\varepsilon,M )r_{F}( d^{\prime},\varepsilon,L ),
$$
and
$$s_{F}(d\times d^{\prime},\varepsilon,M \times L) \geq s_{F}(d,\varepsilon,M )s_{F}(d^{\prime},\varepsilon,L ).
$$
\end{lem}
\begin{proof}
    Fix $\varepsilon>0$ and $F \in \text{Fin}(G)$, let $\{x_{1},x_{2},\cdots,x_{k}\}$ be the $(F,\varepsilon)$- spanning set of $M$ with the smallest cardinality and $\{y_{1},y_{2},\cdots,y_{w}\}$ be the $(F,\varepsilon)$- spanning set of $L$ with the smallest cardinality. Observe that for each $(x,y) \in M \times L$, the point $x \in B_{d_{F}}(x_{i},\varepsilon)$ for some $i \in \{1,2\,\cdots,k\}$ and $y \in B_{d_{F}}(y_{j},\varepsilon)$ for some $j \in \{1,2\,\cdots,w\}$, then
    \begin{equation*}(d\times d^{\prime})_{F}((x,y),(x_{i},y_{j})) =\max\{d_{F}(x,x_{i}) , d^{\prime}_{F}(y,y_{j})\} \leq \varepsilon,
    \end{equation*}
    which implies \begin{equation*}r_{F}(d\times d^{\prime},\varepsilon,M \times L) \leq kw=r_{F}(d,\varepsilon,M )r_{F}( d^{\prime},\varepsilon,L ).\end{equation*}

\hspace{4mm}
Now we prove the second inequality. Let $\{x_{1},x_{2},\cdots,x_{k^{\prime}}\}$ be the $(F,\varepsilon)$- separated set of $M$ with the maximal cardinality and $\{y_{1},y_{2},\cdots,y_{w^{\prime}}\}$ be the $(F,\varepsilon)$- separated set of $L$ with the maximal cardinality. Note that $d_{F}(x_{i},x_{t}) \leq \varepsilon$ implies that $x_{i}=x_{t}$, and $d_{F}(y_{j},y_{h}) \leq \varepsilon$ implies that $y_{j}=y_{h}$. Let $(x_{i},y_{j})$ and $(x_{t},y_{h})$ be distinct elements of the set
\begin{equation*}P=\{(x_{i},y_{j})|i=1,2, \cdots, k^{\prime},~ j =1,2, \cdots, w^{\prime}\} \sub M \times L
\end{equation*}
in which case either $d_{F}(x_{i},x_{t}) > \varepsilon$ or $d_{F}(y_{j},y_{h}) >\varepsilon$ holds. Hence
\begin{equation}\label{wooo}
    (d\times d^{\prime})_{F}((x_{i},y_{j}),(x_{t},y_{h})) =\max\{d_{F}(x_{i},x_{t}) , d^{\prime}_{F}(y_{j},y_{h})\} >\varepsilon.
\end{equation}
In particular, the two balls $B_{(d\times d^{\prime})_{F}}((x_{i},y_{j}),\varepsilon/2)$ and $B_{(d\times d^{\prime})_{F}}((x_{t},y_{h}),\varepsilon/2)$ are disjoint, otherwise  there exists a point $z\in X\times Y$ lies in both balls, then
\begin{equation*}(d\times d^{\prime})_{F}((x_{i},y_{j}),(x_{t},y_{h})) \leq (d\times d^{\prime})_{F}((x_{i},y_{j}),z)+(d\times d^{\prime})_{F}(z,(x_{t},y_{h})) \leq \varepsilon
\end{equation*}
contradicting (\ref{wooo}). 

\hspace{4mm}
Therefore we conclude that $P$ is a $(F,\varepsilon)$- separated set of $M \times L$, then \begin{equation*}s_{F}(d\times d^{\prime},\varepsilon,M \times L) \geq k^{\prime} w^{\prime}=s_{F}(d,\varepsilon,M )s_{F}(d^{\prime},\varepsilon,L ).\qedhere
\end{equation*}
\end{proof}

\begin{lem}\label{2.22}
       Let $(X,G)$ be a $G$-system with a metric $d$. If $M \sub X$ is a compact set, then for any F$\phi$lner sequence $\{F_{n}\}$ in $G$, we have
\begin{align*}
    \underline{\emph{mdim}}_{\emph{M}}(M,G,d)&= \liminf_{\varepsilon \rightarrow 0 }\frac{1}{|\log \varepsilon|}\left(\liminf_{n \rightarrow \infty}\frac{1}{|F_{n}|}r_{F_{n}}(d,\varepsilon,M)\right)\\
    &=\liminf_{\varepsilon \rightarrow 0 }\frac{1}{|\log \varepsilon|}\left(\liminf_{n \rightarrow \infty}\frac{1}{|F_{n}|}s_{F_{n}}(d,\varepsilon,M)\right),
\end{align*}
\begin{align*}
\overline{\emph{mdim}}_{\emph{M}}(M,G,d)&= \limsup_{\varepsilon \rightarrow 0 }\frac{1}{|\log \varepsilon|}\left(\liminf_{n \rightarrow \infty}\frac{1}{|F_{n}|}r_{F_{n}}(d,\varepsilon,M)\right)\\
&=\limsup_{\varepsilon \rightarrow 0 }\frac{1}{|\log \varepsilon|}\left(\liminf_{n \rightarrow \infty}\frac{1}{|F_{n}|}s_{F_{n}}(d,\varepsilon,M)\right).
\end{align*}

\end{lem}
\begin{proof}
We only prove the formula for the lower metric mean dimension. Let $\{F_{n}\}$ be a F$\phi$lner sequence in $G$, according to Remark \ref{01010}, we have
\begin{equation}\label{1000001}
  r_{F_{n}}(d,\varepsilon,M) \leq s_{F_{n}}(d,\varepsilon,M) \leq 
\text{cov}_{F_{n}}(d,\varepsilon,M).  
\end{equation}

\hspace{4mm}
Notice that if $W$ is an $(F_{n},\varepsilon)$-spanning set of $M$, then the $d_{F_{n}}$-balls of radius $\varepsilon$ cover $E$ and these are $|W|$ sets of $d_{F_{n}}$ diameter smaller than $2\varepsilon$. This implies that 
\begin{equation}\label{999099}
    r_{F_{n}}(d,\varepsilon,M) \geq 
\text{cov}_{F_{n}}(d,2\varepsilon,M) .
\end{equation}

\hspace{4mm}
Combining (\ref{1000001}) and (\ref{999099}), we have
\begin{equation*}
\text{cov}_{F_{n}}(d,2\varepsilon,M) \leq r_{F_{n}}(d,\varepsilon,M) \leq s_{F_{n}}(d,\varepsilon,M) \leq 
\text{cov}_{F_{n}}(d,\varepsilon,M)
\end{equation*}
Therefore
\begin{equation*}
\frac{1}{|F_{n}|}\text{cov}_{F_{n}}(d,2\varepsilon,M) \leq \frac{1}{|F_{n}|}r_{F_{n}}(d,\varepsilon,M) \leq \frac{1}{|F_{n}|}s_{F_{n}}(d,\varepsilon,M) \leq 
\frac{1}{|F_{n}|}\text{cov}_{F_{n}}(d,\varepsilon,M).
\end{equation*}
 Taking  the  limit infimum as   $n \to \infty$, combining with Remark \ref{01010}, we have
 \begin{align*}
\lim_{n \rightarrow \infty}\frac{1}{|F_{n}|}\text{cov}_{F_{n}}(d,2\varepsilon,M) &\leq \liminf_{n \rightarrow \infty}\frac{1}{|F_{n}|}r_{F_{n}}(d,\varepsilon,M)\\
&\leq \liminf_{n \rightarrow \infty}\frac{1}{|F_{n}|}s_{F_{n}}(d,\varepsilon,M)\\
&\leq 
\lim_{n \rightarrow \infty}\frac{1}{|F_{n}|}\text{cov}_{F_{n}}(d,\varepsilon,M). 
 \end{align*}
 Hence
  \begin{align*}
\frac{1}{|\log \varepsilon|}\lim_{n \rightarrow \infty}\frac{1}{|F_{n}|}\text{cov}_{F_{n}}(d,2\varepsilon,M) &\leq \frac{1}{|\log \varepsilon|}\liminf_{n \rightarrow \infty}\frac{1}{|F_{n}|}r_{F_{n}}(d,\varepsilon,M)\\
&\leq \frac{1}{|\log \varepsilon|}\liminf_{n \rightarrow \infty}\frac{1}{|F_{n}|}s_{F_{n}}(d,\varepsilon,M)\\
&\leq 
\frac{1}{|\log \varepsilon|}\lim_{n \rightarrow \infty}\frac{1}{|F_{n}|}\text{cov}_{F_{n}}(d,\varepsilon,M). 
 \end{align*}
 Taking the limit infimum as   $\varepsilon \to 0$, we get the desired result.
\end{proof}

\begin{thm}\label{a2.2}
   Let $(X,G)$ and $(Y,G)$ be two $G$-systems, where $(X,d)$ and $(Y,d')$ are compact metric spaces. If $M \sub X$ and $L \sub Y$ are compact sets, then
\begin{align*}
&\underline{\emph{mdim}}_{\emph{M}}(M,G,d) + \underline{\emph{mdim}}_{\emph{M}}(L, G,d')
\leq \underline{\emph{mdim}}_{\emph{M}}(M \times L, G , d \times d^{\prime}) \\
&\leq \min\{\overline{\emph{mdim}}_{\emph{M}}(M,G,d) + \underline{\emph{mdim}}_{\emph{M}}(L, G,d'),\underline{\emph{mdim}}_{\emph{M}}(M,G,d) + \overline{\emph{mdim}}_{\emph{M}}(L, G,d') \}\\
&\leq \max\{\overline{\emph{mdim}}_{\emph{M}}(M,G,d) + \underline{\emph{mdim}}_{\emph{M}}(L, G,d'),\underline{\emph{mdim}}_{\emph{M}}(M,G,d) + \overline{\emph{mdim}}_{\emph{M}}(L, G,d') \}\\
&\leq \overline{\emph{mdim}}_{\emph{M}}( M\times L, G , d \times d^{\prime})\\
&\leq \overline{\emph{mdim}}_{\emph{M}}(M,G,d) + \overline{\emph{mdim}}_{\emph{M}}(L, G,d').
\end{align*}

\end{thm}
\begin{proof}
 Let $\{F_{n}\}$ be any F$\phi$lner sequence in $G$, by Lemma \ref{2.11}, we have
 $$\log r_{F_{n}}(d\times d^{\prime},\varepsilon,M \times L) \leq \log r_{F_{n}}(d,\varepsilon,M )+ \log r_{F_{n}}( d^{\prime},\varepsilon,L ),
$$
 and
 $$\log s_{F_{n}}(d\times d^{\prime},\varepsilon,M \times L) \geq \log s_{F_{n}}(d,\varepsilon,M )+\log s_{F_{n}}( d^{\prime},\varepsilon,L ).
 $$
 Then
$$\frac{1}{|F_{n}|}\log r_{F_{n}}(d\times d^{\prime},\varepsilon,M \times L) \leq \frac{1}{|F_{n}|}\log r_{F_{n}}(d,\varepsilon,M )+\frac{1}{|F_{n}|}\log r_{F_{n}}( d^{\prime},\varepsilon,L ),
 $$
 and
$$\frac{1}{|F_{n}|}\log s_{F_{n}}(d\times d^{\prime},\varepsilon,M \times L) \geq \frac{1}{|F_{n}|}\log s_{F_{n}}(d,\varepsilon,M )+\frac{1}{|F_{n}|}\log s_{F_{n}}( d^{\prime},\varepsilon,L ).
 $$
Hence, taking  the  limit supremum as   $n \to \infty$, we get
 $$\begin{aligned} 
r(d \times d^{\prime},\varepsilon,\{F_{n}\},M \times L)&=\limsup_{n \rightarrow \infty}\frac{1}{|F_{n}|}\log r_{F_{n}}(d\times d^{\prime},\varepsilon,M \times L)\\
&\leq \limsup_{n \rightarrow \infty}\left(\frac{1}{|F_{n}|}r_{F_{n}}(d,\varepsilon,M)+\frac{1}{|F_{n}|}r_{F_{n}}(d^{\prime},\varepsilon,L)\right)\\
&\leq \limsup_{n \rightarrow \infty}\frac{1}{|F_{n}|}r_{F_{n}}(d,\varepsilon,M)+\limsup_{n \rightarrow \infty}\frac{1}{|F_{n}|}r_{F_{n}}(d^{\prime},\varepsilon,L)\\
&=r(d,\varepsilon,\{F_{n}\},M)+r(d^{\prime},\varepsilon,\{F_{n}\},L),
 \end{aligned} $$
 and 
  $$\begin{aligned} 
s(d \times d^{\prime},\varepsilon,\{F_{n}\},M \times L)&=\limsup_{n \rightarrow \infty}\frac{1}{|F_{n}|}\log s_{F_{n}}(d\times d^{\prime},\varepsilon,M \times L)\\
&\geq \limsup_{n \rightarrow \infty}\left(\frac{1}{|F_{n}|}s_{F_{n}}(d,\varepsilon,M)+\frac{1}{|F_{n}|}s_{F_{n}}(d^{\prime},\varepsilon,L)\right)\\
&\geq \liminf_{n \rightarrow \infty}\frac{1}{|F_{n}|}s_{F_{n}}(d,\varepsilon,M)+\limsup_{n \rightarrow \infty}\frac{1}{|F_{n}|}s_{F_{n}}(d^{\prime},\varepsilon,L)\\
&=\liminf_{n \rightarrow \infty}\frac{1}{|F_{n}|}s_{F_{n}}(d,\varepsilon,M)+s(d^{\prime},\varepsilon,\{F_{n}\},L).
 \end{aligned}$$ 
 Therefore
 \begin{align*}
&\overline{\text{mdim}}_{\text{M}}(M,G,d) + \overline{\text{mdim}}_{\text{M}}(L, G,d')\\
   &=\limsup_{\varepsilon \rightarrow 0 }\frac{r(d,\varepsilon,\{F_{n}\},M)}{|\log \varepsilon|}+\limsup_{\varepsilon \rightarrow 0 }\frac{r(d^{\prime},\varepsilon,\{F_{n}\},L)}{|\log \varepsilon|}\\
   &\geq \limsup_{\varepsilon \rightarrow 0 }\left(\frac{r(d,\varepsilon,\{F_{n}\},M)}{|\log \varepsilon|}+\frac{r(d^{\prime},\varepsilon,\{F_{n}\},L)}{|\log \varepsilon|}\right)\\
   &=\limsup_{\varepsilon \rightarrow 0 }\frac{1}{{|\log \varepsilon|}}\left(r(d,\varepsilon,\{F_{n}\},M)+r(d^{\prime},\varepsilon,\{F_{n}\},L)\right)\\
   &\geq \limsup_{\varepsilon \rightarrow 0 }\frac{r(d \times d^{\prime},\varepsilon,\{F_{n}\},M \times L)}{|\log \varepsilon|}\\
&=\overline{\text{mdim}}_{\text{M}}(M \times L, G , d \times d^{\prime}), 
 \end{align*}
 and
\begin{align*}
&\underline{\text{mdim}}_{\text{M}}(M,G,d) + \overline{\text{mdim}}_{\text{M}}(L, G,d')\\
   &=\liminf_{\varepsilon \rightarrow 0 }\frac{r(d,\varepsilon,\{F_{n}\},M)}{|\log \varepsilon|}+\limsup_{\varepsilon \rightarrow 0 }\frac{r(d^{\prime},\varepsilon,\{F_{n}\},L)}{|\log \varepsilon|}\\
   &\geq \liminf_{\varepsilon \rightarrow 0 }\left(\frac{r(d,\varepsilon,\{F_{n}\},M)}{|\log \varepsilon|}+\frac{r(d^{\prime},\varepsilon,\{F_{n}\},L)}{|\log \varepsilon|}\right)\\
   &=\liminf_{\varepsilon \rightarrow 0 }\frac{1}{|\log \varepsilon|}\left(r(d,\varepsilon,\{F_{n}\},M)+r(d^{\prime},\varepsilon,\{F_{n}\},L)\right)\\
   &\geq \liminf_{\varepsilon \rightarrow 0 }\frac{r(d \times d^{\prime},\varepsilon,\{F_{n}\},M \times L)}{|\log \varepsilon|}\\
  & =\underline{\text{mdim}}_{\text{M}}(M \times L, G , d \times d^{\prime}).
\end{align*}

 \hspace{4mm}
On the other hand, 
by Lemma \ref{2.22}, we have
 $$\begin{aligned}
&\overline{\text{mdim}}_{\text{M}}(M \times L, G , d \times d^{\prime})\\
&= \limsup_{\varepsilon \rightarrow 0 }\frac{s(d \times d^{\prime},\varepsilon,\{F_{n}\},M \times L)}{|\log \varepsilon|}\\
   &\geq \limsup_{\varepsilon \rightarrow 0 }\frac{1}{|\log \varepsilon|} \left( \liminf_{n \rightarrow \infty}\frac{1}{|F_{n}|}s_{F_{n}}(d,\varepsilon,M)+s(d^{\prime},\varepsilon,\{F_{n}\},L)\right)\\
   &\geq \liminf_{\varepsilon \rightarrow 0 }\frac{1}{|\log \varepsilon|}\left(\liminf_{n \rightarrow \infty}\frac{1}{|F_{n}|}s_{F_{n}}(d,\varepsilon,M)\right)+\limsup_{\varepsilon \rightarrow 0 }\frac{1}{|\log \varepsilon|}s(d^{\prime},\varepsilon,\{F_{n}\},L)\\
&=\underline{\text{mdim}}_{\text{M}}(M,G,d) + \overline{\text{mdim}}_{\text{M}}(L, G,d'),
 \end{aligned}$$
 and
 \begin{align*}
&\underline{\text{mdim}}_{\text{M}}(M \times L, G , d \times d^{\prime})\\
&= \liminf_{\varepsilon \rightarrow 0 }\frac{s(d \times d^{\prime},\varepsilon,\{F_{n}\},M \times L)}{|\log \varepsilon|}\\
   &\geq \liminf_{\varepsilon \rightarrow 0 }\frac{1}{|\log \varepsilon|} \left( \liminf_{n \rightarrow \infty}\frac{1}{|F_{n}|}s_{F_{n}}(d,\varepsilon,M)+s(d^{\prime},\varepsilon,\{F_{n}\},L)\right)\\
    &\geq \liminf_{\varepsilon \rightarrow 0 }\frac{1}{|\log \varepsilon|}\left(\liminf_{n \rightarrow \infty}\frac{1}{|F_{n}|}s_{F_{n}}(d,\varepsilon,M)\right)+\liminf_{\varepsilon \rightarrow 0 }\frac{1}{|\log \varepsilon|}s(d^{\prime},\varepsilon,\{F_{n}\},L)\\
&=\underline{\text{mdim}}_{\text{M}}(M,G,d) + \underline{\text{mdim}}_{\text{M}}(L, G,d').   \qedhere
 \end{align*}
\end{proof}

\hspace{4mm}
With Theorem \ref{a2.2}, we can easily obtain the following result.
\begin{cor}
    Let $(X,G)$ and $(Y,G)$ be two $G$-systems, where $(X,d)$ and $(Y,d')$ are compact metric spaces. For any compact set $M \sub X$ and $L \sub Y$, if $\underline{\emph{mdim}}_{\emph{M}}(M,G,d)=\overline{\emph{mdim}}_{\emph{M}}(M,G,d)$ and $\underline{\emph{mdim}}_{\emph{M}}(L, G,d')=\overline{\emph{mdim}}_{\emph{M}}(L, G,d')$,
    then
    $$
\emph{mdim}_{\emph{M}}( M\times L, G , d \times d^{\prime})=\emph{mdim}_{\emph{M}}(M,G,d) + \emph{mdim}_{\emph{M}}(L, G,d').$$
\end{cor}

\hspace{4mm}
As an application of the product formula for metric mean dimension, we consider the following example.
\begin{example}
    Let $G$ be a countable discrete amenable group. Suppose $\mathbb{R} / \mathbb{Z}$  is a circle with a metric  $\rho$  defined by
		$$\rho(x, y)=\min _{n \in \mathbb{Z}}|x-y-n|,$$
		and $(\mathbb{R} / \mathbb{Z})^{G}$ is the infinite dimensional torus,
		a metric  $d$ of $(\mathbb{R} / \mathbb{Z})^{G}$ is given by
		$$d\left(\left(x_{g}\right)_{g \in G},\left(y_{g}\right)_{g \in G}\right)=\sum_{g\in G} \alpha_{g} \rho\left(x_{g}, y_{g}\right) ,$$
where
$\alpha_{g} \in (0,+\infty)$ satisfies $$
\alpha_{1_{G}}=1,\sum_{g \in G}\alpha_{g} < +\infty.
$$
  
		The shift map $ \sigma: G\times (\mathbb{R} / \mathbb{Z})^{G} \rightarrow(\mathbb{R} / \mathbb{Z})^{G} $ is defined by 
		$$\sigma^h\left(\left(x_{g}\right)_{g \in G}\right)=\left(x_{gh}\right)_{g \in G}, \text{ for any } h\in G.$$
		For any closed subset $E \subset(\mathbb{R} / \mathbb{Z})^{G}$ satisfying $ \sigma^{h}(E) \subset E$ for any $h\in G$, let $\{F_{n}\}$ be a F$\phi$lner sequence in $G$, then Li and Luo \cite{LL} proved that
	$$\underline{\mathrm{mdim}}_{\mathrm{H}}(E, \{F_{n}\}, d)=\overline{\mathrm{mdim}}_{\mathrm{H}}(E,\{F_{n}\},d) =\underline{\mathrm{mdim}}_{\mathrm{M}}(E,G,d)=\overline{\mathrm{mdim}}_{\mathrm{M}}(E,G,d).$$
	\hspace{4mm}
 According to Theorem \ref{a2.2}, for two closed subsets $M \subset(\mathbb{R} / \mathbb{Z})^{G}$ with $ \sigma^{h}(M) \subset M$ for any $h\in G$, $L \subset(\mathbb{R} / \mathbb{Z})^{G}$ with $ \sigma^{h}(L) \subset L$ any $h\in G$, we have 
 \begin{align*}
\emph{mdim}_{\emph{M}}( M\times L, G , d \times d)
&=\emph{mdim}_{\emph{M}}(M,G,d) + \emph{mdim}_{\emph{M}}(L, G,d)\\
&=\emph{mdim}_{\emph{H}}(M,\{F_{n}\},d) + \emph{mdim}_{\emph{M}}(L, G,d)\\
&=\emph{mdim}_{\emph{M}}(M,G,d) + \emph{mdim}_{\emph{H}}(L, \{F_{n}\},d)
 \end{align*}
 for any F$\phi$lner sequence $\{F_{n}\}$ in $G$.
\end{example}

\hspace{4mm}
This result also reveals the relationship between the product formula for the metric mean dimension and the mean Hausdorff dimension.

\hspace{4mm}
 Next, we prove the product formula for the mean Hausdorff dimension for amenable group actions. We need the following two lemmas.
\begin{lem}\cite{ABB}\label{2001}
    Let  $\varepsilon>0$. Let $(X,d)$ be a compact metric space and $E \sub X$ be a compact set. Suppose   there is a Borel measure $\mu$ on $E$ such that $\mu(E) =1$ and for any open  ball  $E_{i}$   with $\emph{diam}_d (E_{i}) \leq \varepsilon$, we have that
$$
    \mu(E_{i}) \leq (\emph{diam}_d (E_{i}))^{{s}} \quad\text{for any }i\geq 1.
$$
Then, 
$$\emph{dim}^{\star}_{\emph{H}}(E,d,\varepsilon) \geq  {s}.$$
\end{lem}

\begin{lem}\cite{LT}\label{2000}
    Let $c\in (0,1)$. Let $(X,d)$ be a compact metric space and $E \sub X$ be a compact set. There exists $ \varepsilon_0 =\varepsilon_0(c)\in(0,1)$ depending only on $c$ such that for any $0 < \varepsilon \leq \varepsilon_0$, there exists a Borel probability measure $\mu$ on $E$ satisfies
$$\mu(E^{\prime}) \leq (\emph{diam}_d (E^{\prime}))^{c\cdot  \emph{dim}_\emph{H}(E,d,\varepsilon)} $$
for all $E^{\prime} \subset E$ with $\emph{diam}_d (E^{\prime}) < \frac{\varepsilon}{6}.$
\end{lem}

\begin{thm} \label{t38} 
Let $(X,G)$ and $(Y,G)$ be two $G$-systems, where $(X,d)$ and $(Y,d')$ are compact metric spaces. Let $\{F_{n}\}$ be any F$\phi$lner sequence in $G$, if $M \sub X$ and $L \sub Y$ are compact sets, then
$$\underline{\emph{mdim}}_{\emph{H}}(M \times L, \{F_{n} \}, d \times d^{\prime}) \geq \underline{\emph{mdim}}_{\emph{H}}(M,\{F_{n}\},d) + \underline{\emph{mdim}}_{\emph{H}}(L, \{F_{n}\},d'),$$
$$\overline{\emph{mdim}}_{\emph{H}}(M \times L, \{F_{n} \}, d \times d^{\prime}) \geq \overline{\emph{mdim}}_{\emph{H}}(M,\{F_{n}\},d) + \overline{\emph{mdim}}_{\emph{H}}(L, \{F_{n}\},d').$$
\end{thm}

\begin{proof}Fix $0< c < 1$. Let $\{F_{n}\}$ be any F$\phi$lner sequence in $G$.
  It  follows  from Lemma \ref{2000} that there exists $\varepsilon_0 = \varepsilon_0(c) \in (0,1)$ such that for all $0< \varepsilon \leq \varepsilon_{0}$, there are Borel probability measures $\mu$ and $\nu$ in $(M,d)$ and $(L,d')$, respectively, satisfying 
$$\mu(M^{\prime}) \leq (\text{diam}_d(M^{\prime}))^{c\cdot  \text{dim}_{\text{H}}(M,d,\varepsilon)} ,\quad \nu(L^{\prime}) \leq (\text{diam}_{d^{\prime}}(L^{\prime}))^{c\cdot  \text{dim}_{\text{H}}(L,d^{\prime},\varepsilon)}$$
for all $M^{\prime} \subset M$ and $L^{\prime} \subset L$ with $\text{diam}_d (M^{\prime}) < \frac{\varepsilon}{6}$ and  $\text{diam}_{d^{\prime}} (L^{\prime}) < \frac{\varepsilon}{6}$. 

\hspace{4mm}
It is not difficult to verify that $\mu \times \nu$ is the product measure on $M \times L$. Consider the ball $M^{\prime} \times L^{\prime} \sub M \times L$, where $M^{\prime} \subseteq M$ and $L^{\prime} \subseteq L$. Observe that
$$\text{diam}_{d\times d^{\prime}}(M^{\prime} \times L^{\prime}) \geq \max(\text{diam}_d (M^{\prime}), \text{diam}_{d^{\prime}}(L^{\prime})).$$
Then for all $M^{\prime} \times L^{\prime} \sub M \times L $ satisfying $\text{diam}_{d \times d'}(M^{\prime} \times L^{\prime}) < \frac{\varepsilon}{6}$, we have
    \begin{align*}
    (\mu \times \nu) (M^{\prime} \times L^{\prime})
    &= \mu(M^{\prime}) \nu(L^{\prime})\\  &\leq (\text{diam}_d(M^{\prime}))^{c\cdot \text{dim}_{\text{H}}(M,d,\varepsilon)}(\text{diam}_{d'}(L^{\prime}))^{c\cdot \text{dim}_{\text{H}}(L,d',\varepsilon)}\\
    & \leq (\text{diam}_{d\times d'}(M^{\prime} \times L^{\prime}))^{c\cdot  \text{dim}_{\text{H}}(M,d,\varepsilon)}(\text{diam}_{d\times d'}(M \times L^{\prime}))^{c\cdot  \text{dim}_{\text{H}}(L,d',\varepsilon)}\\
    & = (\text{diam}_{d\times d'}(M \times L))^{c\cdot ( \text{dim}_{\text{H}}(M,d,\varepsilon) + \text{dim}_{\text{H}}(L,d',\varepsilon))}.
    \end{align*}
By Lemma \ref{2001}, we have $$\text{dim}_{\text{H}} ^{\star} (M \times L,d \times d^{\prime}, \frac\varepsilon 6) \geq c\cdot ( \text{dim}_{\text{H}}(M,d,\varepsilon) + \text{dim}_{\text{H}}(L,d^{\prime},\varepsilon)).$$ 

\hspace{4mm}
Next, for every $k\geq 1$,  we can take a $ c_{k} \in (0,1)$ such that $c_{k}\rightarrow 1$ as $k\rightarrow \infty$. It follows from the above fact that there exists a $\varepsilon_k =\varepsilon_{k}(c_{k}) \in (0,1)$ such that $\varepsilon_{k}\rightarrow 0$ as $k\rightarrow \infty$ and 
$$\text{dim}_{\text{H}}^{\star}(M \times L,(d \times d^{\prime})_{F_{n}},\frac{ \varepsilon_{k}}{6}) \geq c_{k} ( \text{dim}_{\text{H}}(M,d_{F_{n}},\varepsilon_{k}) + \text{dim}_{\text{H}}(L,d'_{F_{n}},\varepsilon_{k})),$$
for all $ n ,k\in \N$. Hence, we get
$$
\frac{1}{|F_{n}|}\text{dim}_{\text{H}}^{\star}(M \times L,(d \times d')_{F_{n}}, \frac{\varepsilon_{k}}{6}) \geq \frac{c_{k}}{|F_{n}|} ( \text{dim}_{\text{H}}(M,d_{F_{n}},\varepsilon_{k}) + \text{dim}_{\text{H}}(L,d'_{F_{n}},\varepsilon_{k})) .
$$

Therefore,  taking  the  limit infimum and supremum as   $n \to \infty$ and the limit as $k \to \infty$, we have $$\underline{\text{mdim}}_{\text{H}}(M \times L, \{F_{n} \}, d \times d^{\prime}) \geq \underline{\text{mdim}}_{\text{H}}(M,\{F_{n}\},d) + \underline{\text{mdim}}_{\text{H}}(L, \{F_{n}\},d'),$$
$$\overline{\text{mdim}}_{\text{H}}(M \times L
, \{F_{n} \}, d \times d^{\prime}) \geq \overline{\text{mdim}}_{\text{H}}(M,\{F_{n}\},d) + \overline{\text{mdim}}_{\text{H}}(L, \{F_{n}\},d'),$$
which implies the desired result.
\end{proof}
\hspace{4mm}
For a dynamical ststem $(X^{\Z},\sigma)$ where $X$ is a compact metric space and $\sigma$ is the shift map on $X^{\Z}$, given any two points $x=(x_{k})_{k \in \Z},y=({y_{k}})_{k \in \Z} \in X^{\Z}$, consider the metric 
$$\boldsymbol{d}(x,y)=\sum_{k \in \Z}\frac{1}{2^{|k|}}d(x_{k},y_{k}).
$$
 \hspace{4mm}
Recall that the upper and lower Minkowski dimension of $(X,d)$ are defined by 
\begin{align*}
\overline{\text{dim}}_{\text{B}}(X,d)&=\limsup_{\varepsilon\to 0}\frac{\log N(\varepsilon)}{|\log \varepsilon|},\\
\underline{\text{dim}}_{\text{B}}(X,d)&=\liminf_{\varepsilon\to 0}\frac{\log N(\varepsilon)}
{|\log \varepsilon|}.\end{align*}
where $N(\varepsilon)$ denotes the maximal cardinality of an $\varepsilon$-separated set in $(X,d)$ for $\varepsilon>0$. In \cite{VR}, Velozo and Velozo proved that 
$$\overline{\text{mdim}}_{\text{M}}(X^{\Z},\boldsymbol{d},\sigma)=\overline{\text{dim}}_{\text{B}}(X,d),~~\underline{\text{mdim}}_{\text{M}}(X^{\Z},\boldsymbol{d},\sigma)=\underline{\text{dim}}_{\text{B}}(X,d).
$$

  \hspace{4mm}
   Next, we will prove this result for amenable group actions.

\begin{thm}\label{2009}
    Let $(X,G,\sigma)$ be a $G$-system with a metric $d$ satisfying Condition \ref{1000}, the full $G$-shift $\sigma$ on $
		X^{G}$ is   defined by
		$$
  \sigma: G \times X^{G} \rightarrow X^{G}
,(h,(x_{g})_{g \in G}) \mapsto (x_{gh})_{g\in G}. 
		$$
  and the metric $\boldsymbol{d}$ on $X^{G}$ is defined by
		$$\boldsymbol{d}\left((x_g)_{g\in G},(y_g)_{g\in G}\right)=\sum_{g \in G}\alpha_{g}d(x_g,y_g), \text{~for~any~} x=(x_{g})_{g \in G},y=(y_{g})_{g \in G} \in X^{G},$$
where
$\alpha_{g} \in (0,+\infty)$ satisfies $$
\alpha_{1_{G}}=1,\sum_{g \in G}\alpha_{g} < +\infty.
$$
Then we have 
    $$\overline{\emph{mdim}}_{\emph{M}}(X^{G},G,\boldsymbol{d})=\overline{\emph{dim}}_{\emph{B}}(X,d),~~\underline{\emph{mdim}}_{\emph{M}}(X^{G},G,\boldsymbol{d})=\underline{\emph{dim}}_{\emph{B}}(X,d).$$
\end{thm}
\begin{proof}
    We only provide a proof for the first equation, the proof for the second equation follows similarly.

 \hspace{4mm}
 Consider a decreasing sequence $\{\varepsilon_k\}_{k\in\N}$ converging to zero such that $\lim_{k \rightarrow \infty} \frac{\log N(\varepsilon_k)}{|\log \varepsilon_{k}|}= \overline{\text{dim}}_{\text{B}}(X,d)$. Let $P_k=\{p_{1},p_{2},...,p_{h_{k}}\}$ be the maximal collection of $\varepsilon_k$-separated points in $X$. We suppose $\lambda_k$ is the probability measure on $X$ that equidistributes the points in $P_{k}$. Define  
$\mu_{k}=(\lambda_{k})^{\otimes G}$ as the product measure on $X^{G}$.

\hspace{4mm}
Let $\{F_{n}\}$ be a tempered F$\phi$lner sequence in $G$ with 
    $\lim_{n \rightarrow \infty}\frac{|F_{n}|}{\log n}=\infty$. Define
$$A_{({i_{g})_{g \in F_{n}}}}=\{x \in X^{G}: x_{g}=p_{i_{g}} \text{~for~all~} g \in F_{n}\},
$$
where $p_{i_{g}}\in P_{k}$, for all $g \in F_{n}.$

\begin{claim}\label{claim1}
    Let $r <\frac{\varepsilon_{k}}{2}$ and $q \in X^{G}$. Then there exists a unique set $A_{({i_{g})_{g \in F_{n}}}}$ such that
    $$\emph{supp}(\mu_{k})\cap B_{d_{F_{n}}}(q,r) \sub A_{({i_{g})_{g \in F_{n}}}}.
    $$
    \begin{proof}
        Let $x \in \text{supp}(\mu_{k})\cap B_{\boldsymbol{d}_{F_{n}}}(q,r)$, then we have 
        $$d(x_g,q_{g}) \leq \boldsymbol{d}(\sigma^{g}x,\sigma^{g}q) \leq \boldsymbol{d}_{F_{n}}(x,q) \leq r <\frac{\varepsilon_{k}}{2}, ~\forall g \in F_{n}.
        $$
        Since $x \in \text{supp}(\mu_{k})$ we conclude $x_{g} \in P_{k}$. Otherwise the neighbourhood $\mathcal{N}_{g}=\cdots \times X \times U_{g} \times X \times \cdots$ of $x$ has zero $\mu_{k}$-measure, where $U_{g} \sub X$ is an open set containing $x_{g}$ with $U_{g}\cap P_{k}=\emptyset$, which is a contradiction. 

        \hspace{4mm}
        From the choice of $r$ we know that $x_{g}$ can take only one value in $P_{k}$, say $p_{i_{g}}$, where $g \in F_{n}.$
    \end{proof}
\end{claim}
\hspace{4mm}
According to  Claim \ref{claim1}, we get
\begin{equation*}
    \mu_{k}(B_{\boldsymbol{d}_{F_{n}}}(q,r)) \leq \mu_{k}(A_{({i_{g})_{g \in F_{n}}}})=\frac{1}{h_{k}^{|F_{n}|}}=\frac{1}{(N(\varepsilon_{k}))^{|F_{n}|}}
\end{equation*}
for every $q \in X^{G}$ and $ r <\frac{\varepsilon_{k}}{2}$. Then if $A$ satisfies $\mu_{k}(A)>1-\delta$ and $A \sub \bigcup_{i=1}^{L}B_{\boldsymbol{d}_{F_{n}}}(z_{i},r)$, we have
\begin{equation*}
    1-\delta< \mu_{k}(A)\leq \mu_{k}(\bigcup_{i=1}^{L}B_{\boldsymbol{d}_{F_{n}}}(z_{i},r))\leq \frac{L}{(N(\varepsilon_{k}))^{|F_{n}|}}.
\end{equation*}

This implies that $r_{F_{n}}(\boldsymbol{d},\mu,r,\delta,X) \geq (1-\delta)(N(\varepsilon_{k}))^{|F_{n}|}$ for every $\delta \in (0,1)$.
Therefore
$$\limsup_{\varepsilon \rightarrow0}\frac{\sup_{\mu }h_{\mu}^{K}(r,\delta,\{F_{n}\})}{|\log r|} \geq \limsup_{k \rightarrow \infty}\frac{h_{\mu_{k}}^{K}(\varepsilon_{k}/3,\delta,\{F_{n}\})}{|\log (\varepsilon_{k}/3)|} \geq \lim_{k \rightarrow \infty}\frac{\log N(\varepsilon_{k})}{|\log (\varepsilon_{k})|}=\overline{\text{dim}}_{\text{B}}(X,d).
$$By Theorem \ref{1001}, we have
$\overline{\text{mdim}}_{\text{M}}(X,G,\boldsymbol{d})\geq \overline{\text{dim}}_{\text{B}}(X,d).$

\hspace{4mm}
Next, we show that the opposite inequality. Since $X$ is a compact metric space, then $H:=\text{diam}(X) < \infty$. Let $0 <\v < \frac{1}{2}$ and $l=\sum_{g \in G}\alpha_{g} \in (1, +\infty)$, take  $S \in \text{Fin}(G)$ such that 
		$$
		\sum_{g \in G \setminus S}\alpha_{g} \leq \frac{\v}{2H}.
		$$
  \hspace{4mm}
Suppose $B=\{x_i\}_{i=1}^M$ is an $\varepsilon$-separated set of $X$ with the maximum cardinality. Consider the function $f: X\to B,f(x)=x_i$, where $x_i$ is the closest point to $x$ in the subset $B$ (when there are many $x_i$ satisfying the case we take the smallest $i$), obviously $x_i$ is well defined. We can extend $f$ to a measurable function on $X$.  Define the sets $A_i=f^{-1}(x_i)$, and  
$$S_{(i_{t})_{t \in SF_{n}}}=\{y\in X^{G}: y_t\in A_{i_t}\textrm{ for all } t \in SF_{n}\},$$ where $i_t\in \{1,...,M\}$ for each $t \in SF_{n}$. For any $z=(z_{g})_{g \in G},y=(y_{g})_{g \in G} \in S_{(i_{t})_{t \in SF_{n}}}$ , we have
		\begin{align*}			\boldsymbol{d}_{F_{n}}((z_{g})_{g \in G},(y_{g})_{g \in G}) &=\max_{h \in  F_{n}} \boldsymbol{d}((y_{gh})_{g \in G},(z_{gh})_{g \in G}) \\
			&=\max_{h \in  F_{n} }\sum_{g \in G}\alpha_{g} {d(y_{gh},z_{gh})}\\
   &= \max_{h \in  F_{n}} \{\sum_{g \in S}\alpha_{g}{d(y_{gh},z_{gh})}+\sum_{g \in G \setminus S}\alpha_{g} {d(y_{gh},z_{gh})}\} \\
   &\leq \sum_{g \in G} \alpha_{g}  \cdot (2\varepsilon)+\frac{\varepsilon}{2H}\cdot H \\
&<3l\v.
 \end{align*}

   Since the collection of sets $S_{(i_{t})_{t \in SF_{n}}}$ is an open cover of $X^{G}$, then
   $$s_{F_{n}}(\boldsymbol{d},3l\varepsilon,X)\leq M^{|SF_{n}|}=N(\varepsilon)^{|SF_{n}|}.$$ 
   For any measure $\mu$, we have 
   $$r_{F_{n}}(\boldsymbol{d},\mu,r,\delta,X)\leq r_{F_{n}}(\boldsymbol{d},\varepsilon,X) \leq s_{F_{n}}(\boldsymbol{d},\varepsilon,X).$$
   Therefore, we conclude
\begin{align*}\limsup_{\varepsilon \rightarrow 0}\frac{\sup_{\mu }h_{\mu}^{K}(\varepsilon,\delta,\{F_{n}\})}{|\log \varepsilon|} \leq & \lim_{\varepsilon \to 0}\frac{1}{|\log \varepsilon|}\ \cdot \left(\limsup_{n \rightarrow \infty}\left(\frac{1}{|F_{n}|}\log s_{F_{n}}(\boldsymbol{d},\varepsilon,X) \cdot \frac{|F_{n}|}{|SF_{n}|}\right)\right)
\\ \leq&  \limsup_{\varepsilon \to 0}\frac{\log N(\varepsilon)}{| \log (\varepsilon/3l)|}=\overline{\text{dim}}_{\text{B}}(X,d),
\end{align*}
here we used the fact $\lim_{n \rightarrow \infty}\frac{|F_{n}|}{|SF_{n}|}=1$, if $\{F_{n}\}$ is a F$\phi$lner sequence in $G$.

\hspace{4mm}
Combining with Theorem \ref{1001}, we get 
$\overline{\text{mdim}}_{\text{M}}(X,G,\boldsymbol{d})\leq \overline{\text{dim}}_{\text{B}}(X,d).$
\end{proof}

\hspace{4mm}
To better understand this conclusion, we present an example.

\begin{example}
		Let $G$ be a countable discrete amenable group. If $[0,1]^{G}$ is an infinite dimensional cube, we define the shift map $\sigma:G\times [0,1]^{G}\to [0,1]^{G}$ by
		$$
  \sigma: G \times [0,1]^{G} \rightarrow K^{G}
,(h,(x_{g})_{g \in G}) \mapsto (x_{gh})_{g\in G}, 
		$$
		and the metric $d$ on $K^{G}$ is defined by
		$$\textbf{d}\left((x_g)_{g\in G},(y_g)_{g\in G}\right)=\sum_{g \in G}\alpha_{g}|x_g-y_g|,$$
where
$\alpha_{g} \in (0,+\infty)$ satisfies $$
\alpha_{1_{G}}=1,\sum_{g \in G}\alpha_{g} < +\infty.
$$

In \cite{LL}, Li and Luo  proved that for $\{F_{n}\}$ F$\phi$lner sequence in $G$, one has 
$$\mathrm{mdim}\left([0,1]^{G},G\right)=\mathrm{mdim}_{\mathrm{M}}\left([0,1]^{G},G,\boldsymbol{d}\right)=\mathrm{mdim}_{\mathrm{H}}\left([0,1]^{G},\{{F_{n}}\},\boldsymbol{d}\right)=1.$$
Therefore, we have
$$\mathrm{mdim}_{\mathrm{M}}\left([0,1]^{G},G,\boldsymbol{d}\right)=\mathrm{dim}_{\mathrm{B}}\left([0,1],d\right)=1.$$
\end{example}

\hspace{4mm}
We note that Theorem \ref{t38} holds in a special case, so we give the following conjecture under the general amenable group action.

\begin{conjecture}
     Let $(X,G,\sigma)$ be a $G$-system with a metric $d$,  where the shift map $\sigma:G\times X^{G}\to X^{G}$ and the metric $\boldsymbol{d}$ on $X^{G}$ are defined  in Theorem \ref{2009}. Then we have 
$$\overline{\emph{mdim}}_{\emph{M}}(X^{G},G,\boldsymbol{d})=\overline{\emph{dim}}_{\emph{B}}(X,d),~~\underline{\emph{mdim}}_{\emph{M}}(X^{G},G,\boldsymbol{d})=\underline{\emph{dim}}_{\emph{B}}(X,d).$$ 
\end{conjecture}

\section{On the continuity of metric mean dimension and Hausdorff mean dimension maps}\label{sec4}

\hspace{4mm}
In this section, we will work with a fixed metrizable compact topological space $(X, \tau)$. We define
$$X(\tau) = \lbrace d \colon d \text{ is a metric for } X \text{ and } \tau_d = \tau \rbrace,$$
where $\tau_{d}$ is the topology induced by $d$ on $X$, and  $X(\tau)$  is endowed with the metric  
$$ D(d_1,d_2) = \max_{x,y \in X} \left\{| d_1(x,y) - d_2(x,y)|: \text{ for } d_1,d_2 \in X(\tau) \right\}.$$  We are going to study the continuity of metric mean dimension and mean Hausdorff dimension on $X(\tau)$, respectively.

\hspace{4mm}
We recall that two metrics on a space $X$ are equivalent if they induce the same topology on $X$. Therefore, if $d$ is a fixed metric on $X$ which induces the topology $\tau$, then $X(\tau)$ consists of all the metrics on $X$ which are equivalent to $d$.

\hspace{4mm}
Let $G$ be a countable discrete group and $\{F_{n}\}$ be any F$\phi$lner sequence in $G$. Fix a continuous action $T:G \times X \rightarrow X$,
consider the functions   \begin{equation*} 
\begin{aligned}
\overline{\text{mdim}}_{\text{M}}(X,G,T) \colon X(\tau)& \to \mathbb{R} \cup \lbrace \infty \rbrace\\
d& \mapsto \overline{\text{mdim}}_{\text{M}}(X,G,d),
\end{aligned} 
\end{equation*}
\begin{equation*} 
\begin{aligned}
\underline{\text{mdim}}_{\text{M}}(X,G,T) \colon X(\tau)& \to \mathbb{R} \cup \lbrace \infty \rbrace\\
d& \mapsto \underline{\text{mdim}}_{\text{M}}(X,G,d),
\end{aligned} 
\end{equation*}
and
\begin{equation*}
\begin{aligned}
    \overline{\text{mdim}}_{\text{H}}(X,\{F_{n}\},T) \colon X(\tau) &\to \mathbb{R} \cup \lbrace \infty \rbrace\\
d& \mapsto \overline{\text{mdim}}_{\text{H}}(X,\{F_{n}\},d),
\end{aligned}
\end{equation*}
\begin{equation*}
\begin{aligned}
    \underline{\text{mdim}}_{\text{H}}(X,\{F_{n}\},T) \colon X(\tau) &\to \mathbb{R} \cup \lbrace \infty \rbrace\\
d& \mapsto \underline{\text{mdim}}_{\text{H}}(X,\{F_{n}\},d).
\end{aligned}
\end{equation*}

\begin{remark}\label{3000}
   Recall that two metrics $d_{1}$ and $d_{2}$ on $X$ are \textbf{uniformly equivalent} if there are two real constants $0<a \leq b$ such that
   $$ad_{1}(x,y) \leq d_{2}(x,y) \leq bd_{1}(x,y)
   $$
   for every $x,y \in X$. Let $\{F_{n}\}$ be a F$\phi$lner sequence in $G$. If $d_1$ and $d_2$ on $X$ are two uniformly equivalent metrics on $X$, it is not difficult to see  that
$$\emph{mdim}_{\emph{M}}(X , G,d_1) = \emph{mdim}_{\emph{M}}(X , G,d_2), \quad  
\emph{mdim}_{\emph{H}}(X,\{F_{n}\},d_1) =  \emph{mdim}_{\emph{H}}(X,\{F_{n}\},d_2).$$  
\end{remark}

 \begin{remark} Note that if $h_{\text{top}}(X,G) < \infty$, then $\emph{mdim}_{\emph{M}}(X,G,d)=0$. As the topological entropy does not depend on
the metrics, then $\emph{mdim}_{\emph{M}}(X,G,\tilde{d})=0$ for any $\tilde{d} \in X(\tau$). Moreover, for any F$\phi$lner sequence $\{F_{n}\}$ in $G$, according to Remark \ref{x000}, we have $\emph{mdim}_{\emph{H}}(X,\{F_{n}\},\tilde{d})=0$ for any $\tilde{d} \in X(\tau$). Hence, if $h_{\text{top}}(X,G) < \infty$, then  $$\overline{\emph{mdim}}_{\emph{M}}(X,G,T) \colon X(\tau) \to \mathbb{R} \quad \text{ and }\quad \overline{\emph{mdim}}_{\emph{H}}(X,\{F_{n}\},T) \colon X(\tau) \to \mathbb{R} $$ are  the zero maps. 
  \end{remark}

\hspace{4mm}
In the following theorem, we will consider the continuity of metric mean dimension in $X(\tau)$ for amenable group actions.
 \begin{thm}\label{3002} Let $(X,G)$ be a $G$-system. If there exists a continuous action  $T:G \times X\rightarrow X$ such that $\overline{\emph{mdim}}_{\emph{M}}(X,G,d)>0$ for some $d\in X(\tau )$, then 
\begin{equation*}
\begin{aligned}
\overline{\emph{mdim}}_{\emph{M}}(X,G,T) \colon X(\tau) &\to \mathbb{R} \cup \lbrace \infty \rbrace\\
d &\mapsto \overline{\emph{mdim}}_{\emph{M}}(X,G,d)
\end{aligned}
\end{equation*}
is not continuous anywhere. 
\end{thm}
\begin{proof}   
 Given any $\alpha,\varepsilon\in (0,1)$, we define the metric
\begin{align*}d_{\alpha,\varepsilon}(x,y)=\left\lbrace\begin{array}{cccc}
d(x,y),&\text{ if } d(x,y)\geq \varepsilon,\\
\varepsilon^{1-\alpha}d(x,y)^\alpha,&\text{ if } d(x,y)<\varepsilon.
\end{array}\right.\end{align*}
It's obvious that  $d_{\alpha,\varepsilon}\in X(\tau)$. Moreover, taking $x,y\in X$ such that $d(x,y)\geq \varepsilon$, we have $|d(x,y)-d_{\alpha,\varepsilon}(x,y)|=0<\varepsilon$. Consider $x,y\in X$ satisfying $d(x,y)<\varepsilon$, then
\begin{eqnarray*}
|d(x,y)-d_{\alpha,\varepsilon}(x,y)|=|d(x,y)-\varepsilon^{1-\alpha} d(x,y)^\alpha |
\leq d(x,y) + \varepsilon^{1-\alpha} d(x,y)^\alpha < 2\varepsilon.
\end{eqnarray*}

Therefore, $D(d,d_{\alpha,\varepsilon})<2\varepsilon$.

\hspace{4mm}
Next, we will prove that
$$\overline{\text{mdim}}_{\text{M}}(X,G,d_{\alpha,\varepsilon})=\frac{\overline{\text{mdim}}_{\text{M}}(X,G,d)}{\alpha}.$$
Let $\{F_{n}\}$ be a F$\phi$lner sequence in $G$. Consider any $\eta\in(0,\varepsilon)$. Let $A$ an $(F_{n},\eta)$- spanning set of $(X,d)$. Then for any $y\in X$, there is $x\in A$ satisfying $d_{F_{n}}(x,y)<\eta$. Hence,
 $$(d_{\alpha,\varepsilon})_{F_{n}}(x,y)=\varepsilon^{1-\alpha}d_{F_{n}}(x,y)^\alpha<\varepsilon^{1-\alpha}\eta^\alpha,$$ 
which implies that $A$ is an $(F_{n},\varepsilon^{1-\alpha}\eta^\alpha)$-spanning set of $(X,d_{\alpha,\varepsilon})$. It follows that 
 $$r_{F_{n}}(d_{\alpha,\varepsilon},\varepsilon^{1-\alpha}\eta^\alpha,X)\leq r_{F_{n}}(d,\eta,X).$$
Therefore, we obtain that
\begin{eqnarray}\label{1a}
\overline{\text{mdim}}_{\text{M}}(X,G,d_{\alpha,\varepsilon})&=&\displaystyle\limsup_{\eta\to 0}\limsup_{n \to \infty} \frac{r_{F_{n}}(d_{\alpha,\varepsilon},\varepsilon^{1-\alpha}\eta^\alpha ,X)}{|F_{n}||\log(\varepsilon^{1-\alpha} \eta^\alpha) |} \nonumber\\
&\leq & \displaystyle\limsup_{\eta\to 0}\limsup_{n \to \infty} \frac{r_{F_{n}}(d,\eta,X)}{\alpha|F_{n}||\log \eta |} \frac{|\log (\eta^\alpha)|}{|\log (\varepsilon^{1-\alpha}\eta^\alpha)|}\nonumber\\
&=& \frac{\overline{\text{mdim}}_{\text{M}}(X,G,
d)}{\alpha}.
\end{eqnarray}

\hspace{4mm}
On the other hand, for any $x,y\in X$ satisfying   $(d_{\alpha,\varepsilon})_{F_{n}}(x,y)<\varepsilon$, we have that $d_{F_{n}}(x,y)<\varepsilon$. Let $E$ be an $(F_{n},\eta)$-spanning set of $(X,d_{\alpha,\varepsilon})$, where $\eta\in (0,\varepsilon)$.  Then for any $y\in X$, there is $x\in E$ with $(d_{\alpha,\varepsilon})_{F_{n}}(x,y)<\eta$ and it follows that
$$(d_{\alpha,\varepsilon})_{F_{n}}(x,y)=\varepsilon^{1-\alpha}d_{F_{n}}(x,y)^\alpha<\eta<\varepsilon \Rightarrow  d_{F_{n}}(x,y)<\varepsilon^{\frac{\alpha-1}{\alpha}}\eta^{\frac{1}{\alpha}},$$
which implies that $E$ is an $(F_{n},\varepsilon^{\frac{\alpha-1}{\alpha}}\eta^{\frac{1}{\alpha}})$ spanning set of $(X,d)$. Hence,
$$r_{F_{n}}(d_{\alpha,\varepsilon},\eta,X)\geq r_{F_{n}}(d,\varepsilon^{\frac{1-\alpha}{\alpha}}\eta^{\frac{1}{\alpha}},X).$$
Then we obtain that 
\begin{eqnarray}\label{1ab}
\overline{\text{mdim}}_{\text{M}}(X,G,d_{\alpha,\varepsilon})&=&\displaystyle\limsup_{\eta\to 0}\limsup_{n \to \infty} \frac{r_{F_{n}}(d_{\alpha,\varepsilon},\eta,X)}{|F_{n}||\log \eta |} \nonumber\\
&\geq & \displaystyle\limsup_{\eta\to 0}\limsup_{n \to \infty} \frac{r_{F_{n}}(d,\varepsilon^{\frac{1-\alpha}{\alpha}}\eta^{\frac{1}{\alpha}},X)}{|F_{n}||\log (\varepsilon^{\frac{\alpha-1}{\alpha}}\eta^{\frac{1}{\alpha}}) |}\frac{|\log (\varepsilon^{\frac{\alpha-1}{\alpha}}\eta^{\frac{1}{\alpha}}) |}{|\log \eta|}\nonumber\\
&=&\displaystyle\limsup_{\eta\to 0}\limsup_{n \to \infty} \frac{r_{F_{n}}(d,\varepsilon^{\frac{1-\alpha}{\alpha}}\eta^{\frac{1}{\alpha}},X)}{|F_{n}||\log (\varepsilon^{\frac{\alpha-1}{\alpha}}\eta^{\frac{1}{\alpha}}) |}\frac{|\log (\eta^{\frac{1}{\alpha}}) |}{|\log \eta|}\nonumber\\
&=& \frac{\overline{\text{mdim}}_{\text{M}}(X,G,d)}{\alpha}.
\end{eqnarray}
Combining \eqref{1a} and \eqref{1ab} yields
$\overline{\text{mdim}}_{\text{M}}(X,G,d_{\alpha,\varepsilon})=\frac{\overline{\text{mdim}}_{\text{M}}(X,G,d)}{\alpha}. $

\hspace{4mm}
Given that 
$$\overline{\text{mdim}}_{\text{M}}(X ,G,d_{\alpha,\varepsilon})= \frac{\overline{\text{mdim}}_{\text{M}}(X,G,d)}{\alpha},$$ 
and $D(d_{\alpha,\varepsilon},d)<2\varepsilon$, for any $\varepsilon>0$, we conclude that $\overline{\text{mdim}}_{\text{M}}(X,G,d)$ is not continuous with respect to the metric.
\end{proof}

\hspace{4mm}
Using a similar method, we can prove the following result. 
\begin{thm}
Let $(X,G)$ be a $G$-system. If there exists a continuous action  $T:G \times X\rightarrow X$ such that if $\underline{\emph{mdim}}_{\emph{M}}(X,G,d)>0$ for some $d\in X(\tau )$, then 
\begin{equation*}
\begin{aligned}
\underline{\emph{mdim}}_{\emph{M}}(X,G,T) \colon X(\tau) &\to \mathbb{R} \cup \lbrace \infty \rbrace\\
d &\mapsto \underline{\emph{mdim}}_{\emph{M}}(X,G,d)
\end{aligned}
\end{equation*}
is not continuous anywhere.\end{thm} 

\hspace{4mm}
Next, we consider the continuity of the mean Hausdorff dimension in $X(\tau)$ for amenable group actions. We need the following lemma.
\begin{lem}\label{3003}
    Let $(X,G)$ be a $G$-system with a metric $d$. Fix any $a \in (0, 1]$. Consider the function $\zeta(x) = x^a, x \in [0,\infty)$. Let $\zeta_d(x,y) = \zeta (d(x,y))$, then for any F$\phi$lner sequence $\{F_{n}\}$ in $G$, we have 
\begin{align*}
\overline{\emph{mdim}}_{\emph{H}}(X ,\{F_{n}\},\zeta_{d})&=\frac{\overline{\emph{mdim}}_{\emph{H}}(X,\{F_{n}\},d)}{a},\\ \underline{\emph{mdim}}_{\emph{H}}(X ,\{F_{n}\},\zeta_{d})&=\frac{\underline{\emph{mdim}}_{\emph{H}}(X,\{F_{n}\},d)}{a}.
\end{align*}
\end{lem}
\begin{proof}
 Fix an $a \in (0,1]$, it's clear that $\zeta_d(x,y) = \zeta (d(x,y))$ is a metric on $X$. Given any $\eta>0$, we have that  $d(x,y)\leq \eta$ if and only if  $d(x,y)^a\leq \eta^a$. Hence, it follows that
\begin{align*}
   &\text{H}_{\eta^a}^{s} (X,(\zeta_{d})_{F_{n}})\\
   &= \inf\left\{ \Sigma_{k=1}^{\infty}(\text{diam}_{d_{F_{n}}^a} (E_{k}))^{s}: X=\cup_{k=1}^{\infty} E_{k} \text{ with } \text{diam}_{d_{F_{n}}^a} (E_{k})<\eta^a\text{ for all }k\geq 1\right\}\\
&=  \inf\left\{ \Sigma_{k=1}^{\infty}(\text{diam}_{d_{F_{n}}^a} (E_{k}))^{s}: X=\cup_{k=1}^{\infty} E_{k} \text{ with } \text{diam}_{d_{F_{n}}} (E_{k})<\eta\text{ for all }k\geq 1\right\}\\
&=\inf\left\{ \Sigma_{k=1}^{\infty}(\text{diam}_{d_{F_{n}}} (E_{k}))^{as}: X=\cup_{k=1}^{\infty} E_{k} \text{ with } \text{diam}_{d_{F_{n}}} (E_{k})<\eta\text{ for all }k\geq 1\right\}\\
&=\text{H}_{\eta}^{as} (X,d_{F_{n}}). 
\end{align*}

Therefore,
\begin{align*}
  \dim_{\text{H}}(X,(\zeta_{d})_{F_{n}},\eta^a)&=\sup\{s\geq 0:\text{H}_{\eta^a}^{s} (X,(\zeta_{d})_{F_{n}})\geq 1\}= \sup\{s\geq 0:\text{H}_{\eta}^{a s} (X,d_{F_{n}})\geq 1\}\\
&= \frac{1}{a}\sup\{a s\geq 0:\text{H}_{\eta}^{a s} (X,d_{F_{n}})\geq 1\} = \frac{1}{a} \dim_{\text{H}}(X,d_{F_{n}},\eta). 
\end{align*}

\hspace{4mm}
By the definition of the upper and lower mean Hausdorff dimension, we get the desired result.
\end{proof}

\begin{thm}\label{3004} Let $(X,G)$ be a $G$-system and $\{F_{n}\}$ be a  F$\phi$lner sequence in $G$. If there exists a continuous action  $T:G \times X\rightarrow X$ such that $\overline{\emph{mdim}}_{\emph{H}}(X,\{F_{n}\},d)>0$ for some $d\in X(\tau )$, then 
\begin{equation*}
\begin{aligned}
\overline{\emph{mdim}}_{\emph{H}}(X,\{F_{n}\},T) \colon X(\tau) &\to \mathbb{R} \cup \lbrace \infty \rbrace\\
d &\mapsto \overline{\emph{mdim}}_{\emph{H}}(X,\{F_{n}\},d)
\end{aligned}
\end{equation*}
is not continuous anywhere. 
\end{thm}

\begin{proof}
Let $\{F_{n}\}$ be any F$\phi$lner sequence in $G$. Given any $\alpha,\varepsilon\in (0,1)$, we define the metric
$$d_{\alpha,\varepsilon}(x,y)=\left\lbrace\begin{array}{cccc}
d(x,y),&\text{if } d(x,y)\geq \varepsilon,\\
\varepsilon^{1-\alpha}d(x,y)^\alpha,&\text{if } d(x,y)<\varepsilon.
\end{array}\right.$$

According to the proof of Theorem \ref{3002}, we have $D(d,d_{\alpha,\varepsilon})<2\varepsilon$. 

\hspace{4mm}
By Lemma \ref{3003}, we have the relation $$\overline{\text{mdim}}_{\text{H}}(X,\{F_{n}\},d^{\alpha})=\frac{\overline{\text{mdim}}_{\text{H}}(X,\{F_{n}\},d)}{\alpha},\quad\text{ for any }\alpha\in(0,1).$$

Fix $\eta\in (0,\varepsilon)$. For each $x,y\in X$ with $d_{F_{n}}(x,y)<\eta,$  we have $$(d_{\alpha,\epsilon})_{{F_{n}}}(x,y)=\varepsilon^{1-\alpha}d_{F_{n}}(x,y)^\alpha.$$

Therefore, for all $E\subset M$ such that $\text{diam}_{d_{F_{n}}^\alpha}(E)<\eta$, we have $\text{diam}_{(d_{\alpha,\varepsilon})_{F_{n}}}(E)<\varepsilon^{1-\alpha}\eta$ which implies that
$$\text{H}_{\varepsilon^{1-\alpha}\eta}^s(X,(d_{\alpha,\varepsilon})_{{F_{n}}})\leq \text{H}_\eta^s(X,d_{F_{n}}^\alpha), \quad\text{for every }0<\eta<\epsilon.$$ 
Therefore,
\begin{equation}\label{ncdmfezx}\overline{\text{mdim}}_{\text{H}}(X,\{F_{n}\},d_{\alpha,\varepsilon})\leq \overline{\text{mdim}}_{\text{H}}(X,\{F_{n}\},d^\alpha ) = \frac{\overline{\text{mdim}}_{\text{H}}(X,\{F_{n}\},d)}{\alpha}.\end{equation}  

\hspace{4mm}
On the other hand, given $\eta\in (0,\varepsilon)$,  for each $x,y\in X$ with  $d_{{F_{n}}}(x,y)<\eta$,  we have   
$$(d_{\alpha,\epsilon})_{{F_{n}}}(x,y)= \varepsilon^{1-\alpha} d_{F_{n}}(x,y)^\alpha> \eta^{1-\alpha} d_{F_{n}}(x,y)^\alpha.
 $$ Therefore, for all $E\sub X$ with $\text{diam}_{(d_{\alpha,\varepsilon})_{{F_{n}}}}(E)<\eta$, it follows that $\text{diam}_{d_{F_{n}}^\alpha}(E)<\eta^\alpha$ which implies that $$\text{H}_\eta^s(X,(d_{\alpha,\varepsilon})_{{F_{n}}})\geq \text{H}_{\eta^\alpha}^s(X,d_{F_{n}}^\alpha).$$ 
Hence,
\begin{equation}\label{zxxxm}\overline{\text{mdim}}_{\text{H}}(X ,\{F_{n}\},d_{\alpha,\varepsilon})\geq\overline{\text{mdim}}_{\text{H}}(X ,\{F_{n}\},d^\alpha) = \frac{\overline{\text{mdim}}_{\text{H}}(X,\{F_{n}\},d)}{\alpha}.\end{equation}
It follows from \eqref{ncdmfezx} and \eqref{zxxxm} that $\overline{\text{mdim}}_{\text{H}}(X ,\{F_{n}\},d_{\alpha,\varepsilon})=\frac{\overline{\text{mdim}}_{\text{H}}(X,\{F_{n}\},d)}{\alpha}.$

\hspace{4mm}
Given that 
$$ \overline{\text{mdim}}_{\text{H}}(X,\{F_{n}\},d_{\alpha,\varepsilon})= \frac{\overline{\text{mdim}}_{\text{H}}(X,\{F_{n}\},d)}{\alpha},$$ 
and $D(d_{\alpha,\varepsilon},d)<2\varepsilon$, for any $\varepsilon>0$, we conclude that $\overline{\text{mdim}}_{\text{H}}(X,\{F_{n}\},d)$ is not continuous with respect to the metric.
\end{proof}

\hspace{4mm}
Similarly, we can prove the following conclusion.
\begin{thm}
Let $(X,G)$ be a $G$-system and $\{F_{n}\}$ be a  F$\phi$lner sequence in $G$. If there exists a continuous action  $T:G \times X\rightarrow X$ such that  if $\underline{\emph{mdim}}_{\emph{H}}(X,\{F_{n}\},d)>0$, for some $d\in X(\tau )$, then 
\begin{equation*}
\begin{aligned}
\underline{\emph{mdim}}_{\emph{H}}(X,\{F_{n}\},T) \colon X(\tau) &\to \mathbb{R} \cup \lbrace \infty \rbrace\\
d &\mapsto \underline{\emph{mdim}}_{\emph{H}}(X,\{F_{n}\},d)
\end{aligned}
\end{equation*}
is not continuous anywhere. \end{thm}

\section{Composing metrics with subadditive continuous maps} \label{sec5}
\hspace{4mm}
In this section, we will consider the continuity of metric mean dimension concerning metrics in the set  
\begin{equation*}\label{SAM}
\mathcal{A}_{d}(X) = \lbrace \zeta_d \colon \zeta_d(x,y) = \zeta (d(x,y))  \text{ for all }x, y\in  X,  \text{ and } \zeta \in \mathcal{A}[0,\rho]  \rbrace,
\end{equation*}
where $\rho$ is the diameter of $X$ and \begin{equation*}\label{SA}
\mathcal{A}[0, \rho] = \left\{ \zeta :[0,\rho]\rightarrow [0,\infty) : \zeta\text{ is continuous,  increasing, subadditive  and }\zeta^{-1}( 0) =\{ 0\}   \right\}.
\end{equation*}
\begin{remark}
  The function $\zeta:[0,\infty) \rightarrow [0,\infty)$ is called \textbf{subadditive} if  $\zeta(x+y) \leq \zeta(x)+\zeta(y)$ for all $x,y \in [0,\infty)$.
\end{remark}

\begin{lem}\label{metricg}
For any $\zeta \in\mathcal{A}[0,\rho]$, we have  that:   
\begin{itemize}  \item[(i)] $\zeta_{d}$ is a metric on $X$. 
\item[(ii)] $\zeta_d \in X(\tau)$. Consequently, $\mathcal{A}_{d}(X) \subseteq  X(\tau)$.
\item[(iii)] If $(X,G)$ is a $G$-system, for any F$\phi$lner sequence $\{F_{n}\}$ in $G$ and $x,y\in X$, we have $({\zeta_{d}})_{F_{n}}(x,y)= \zeta(d_{F_{n}}(x,y))$.
\end{itemize}
\end{lem}

\begin{proof}   
The statements of (i) and (ii) follow from the proof of Lemma 6.1 in \cite{ABB}, we only prove (iii). 

\hspace{4mm}
Let $\{F_{n}\}$ be a F$\phi$lner sequence in $G$. Since $\zeta$ is increasing, for any $x,y\in X$, we have that 
\begin{align*}
    (\zeta_d)_{F_{n}}(x,y)  &= \max_{ g \in F_{n} } \lbrace \zeta_d(gx,gy)\rbrace\\
&= \max_{ g \in F_{n} } \lbrace \zeta(d(gx,gy)) \rbrace\\
  &=\zeta (\max_{ g \in F_{n} } \lbrace d(gx,gy) \rbrace)  = \zeta({d}_{F_{n}}(x,y)). \qedhere
\end{align*}
\end{proof}

\hspace{4mm}
Next, we consider the continuity of metric mean dimension for amenable group actions with metrics in  $\mathcal{A}_{d}(X)$. For any continuous map $\zeta\in \mathcal{A}[0, \rho]$, take
$$k_m(\zeta) = \liminf_{\varepsilon \to 0^+} \frac{\log(\zeta(\varepsilon))}{\log(\varepsilon)} , k_M(\zeta) = \limsup_{\varepsilon \to 0^+} \frac{\log(\zeta(\varepsilon))}{\log(\varepsilon)}.$$ 

\begin{lem}\label{3999}     For any $\zeta\in \mathcal{A}[0, \rho]$, we have that  $k_m(\zeta) \leq k_M(\zeta)\leq 1$. 
\end{lem}
\begin{proof}
    The proof of this statement follows from \cite{ABB}, we omit it here.
\end{proof}

\begin{pro}\label{4000}
Let $(X,G)$ be a $G$-system with a metric $d$. Taking $\zeta\in \mathcal{A} [0, \rho] $ such that $k_{m}(\zeta),k_{M}(\zeta)>0$. Set  $\zeta_d(x,y) = \zeta \circ d(x, y)$ for every $x,y \in X$, then
\begin{itemize}
\item[(i)] $ \underline{\emph{mdim}}_{\emph{M}}(X,G,d) \geq k_m(\zeta)\underline{\emph{mdim}}_{\emph{M}}(X,G,\zeta_d).$
\item[(ii)] $ \overline{\emph{mdim}}_{\emph{M}}(X,G,d) \leq k_M(\zeta)\overline{\emph{mdim}}_{\emph{M}}(X,G,\zeta_d).$
\end{itemize}
\end{pro}

\begin{proof} By Lemma \ref{3999}, we suppose that $k_{m}(\zeta),k_{M}(\zeta)\in (0,1]$. Let $\{F_{n}\}$ be a F$\phi$lner sequence in $G$.

(i) Fix  $\varepsilon > 0$. If $d_{F_{n}}(x, y) < \varepsilon$, since $\zeta$ is increasing, then we have that $(\zeta_{d})_{F_{n}}(x, y)=\zeta(d_{F_{n}}(x, y)) \leq \zeta(\varepsilon)$.    Hence, we know that any $(F_{n}, \varepsilon)$-spanning subset with respect to $d$ is an $(F_{n}, \zeta(\varepsilon))$-spanning subset with respect to $\zeta_d$. Thus,
\begin{equation}\label{mfef}r_{F_{n}}(d,\varepsilon,X)\geq r_{F_{n}}(\zeta_{d},\zeta(\varepsilon),X).
\end{equation} 
Since $\zeta$ is continuous and $\zeta(0) = 0$, we have $\underset{\varepsilon \to 0}{\lim} \zeta(\varepsilon)  =0$. Therefore,
\begin{align*}
\underline{\text{mdim}}_{\text{M}}(X,G,{d}) & = \liminf_{\varepsilon \to 0} \limsup_{n \to \infty}   \frac{\log r_{F_{n}}(d,\varepsilon,X)}{|F_{n}||\log(\varepsilon)|}\\
& = \liminf_{\varepsilon \to 0} \limsup_{n \to \infty}  \frac{\log r_{F_{n}}(d,\varepsilon,X)}{|F_{n}||\log(\varepsilon)|} \frac{|\log(\zeta(\varepsilon))|}{|\log(\zeta(\varepsilon))|}\\
 & \geq \liminf_{\varepsilon \to 0} \limsup_{n \to \infty}  \frac{\log r_{F_{n}}(\zeta_d,\zeta(\varepsilon),X)}{|F_{n}||\log(\zeta(\varepsilon))|} \frac{|\log(\zeta(\varepsilon))|}{| \log (\varepsilon)|}\quad(\text{ by } \eqref{mfef})\\
\quad & \geq k_{m}(\zeta)  \liminf_{\varepsilon \to 0} \limsup_{n \to \infty}  \frac{\log r_{F_{n}}(\zeta_d,\zeta(\varepsilon),X)}{|F_{n}||\log(\zeta(\varepsilon))|}\\
& = k_{m}(\zeta)  \underline{\text{mdim}}_{\text{M}}(X,G,\zeta_{d}).
\end{align*}

(ii)   Fix  $\varepsilon > 0$. Let $A$ be an $(F_{n},\varepsilon)$-separated subset with respect to $d$, then for any $x,y \in A$ with $x \neq y$, we have 
$
d_{F_{n}}(x,y) =  \underset{g \in F_{n}}{\max} \{ d(gx,gy)\} > \varepsilon
$. 
Therefore, there is $g_{0} \in F_{n}$ such that $ 
d(g_{0}x,g_{0}y) > \varepsilon.$ Since $\zeta$ is increasing, it follows that
$
\zeta\left( d(g_{0}x,g_{0}y)\right) \geq \zeta(\varepsilon)$ which implies that
$$(\zeta_{d})_{F_{n}}(x,y) = \max_{g \in F_{n}}\left\{ \zeta \left(d(gx,gy) \right)\right\} \geq \zeta(\varepsilon).$$
Thus, $A$ is an $(F_{n}, \zeta(\varepsilon))$-separated subset with respect to $\zeta_d$, then
\begin{equation}\label{mbsgdhd}s_{F_{n}}(d,\varepsilon,X) \leq s_{F_{n}}(\zeta_d,\zeta(\varepsilon),X).
\end{equation}
Hence,
\begin{align*}
\overline{\text{mdim}}_{\text{M}}(X,G,d) & = \limsup_{\varepsilon \to 0} \limsup_{n \to \infty}  \frac{\log s_{F_{n}}(d,\varepsilon,X)}{|F_{n}||\log(\varepsilon)|} \\
&= \limsup_{\varepsilon \to 0} \limsup_{n \to \infty}  \frac{\log s_{F_{n}}(d,\varepsilon,X)}{|F_{n}||\log(\varepsilon)|} \frac{|\log(\zeta(\varepsilon))|}{|\log(\zeta(\varepsilon))|}\\
& \leq \limsup_{\varepsilon \to 0} \limsup_{n \to \infty}  \frac{\log s_{F_{n}}(\zeta_d,\zeta(\varepsilon),X)}{|F_{n}||\log(\zeta(\varepsilon))|} \frac{|\log(\zeta(\varepsilon))|}{|\log(\varepsilon)|}\quad(\text{ by } \eqref{mbsgdhd})\quad \\
\quad & \leq k_{M}(\zeta)  \limsup_{\varepsilon \to 0} \limsup_{n \to \infty} \frac{\log s_{F_{n}}(\zeta_d,\zeta(\varepsilon),X)}{|F_{n}||\log(\zeta(\varepsilon))|}  \\
&=k_{M} (\zeta) \overline{\text{mdim}}_{\text{M}}(X,G,\zeta_d).
\end{align*} 
\hspace{4mm}
Therefore, we obtain that $\overline{\text{mdim}}_{\text{M}}(X,G,d) \leq k_{M} (\zeta) \overline{\text{mdim}}_{\text{M}}(X,G,\zeta_d).$
\end{proof}

\begin{lem} \label{4004}
    For any $\zeta \in \mathcal{A}[0, \rho]$ satisfying $k(\zeta) = k_m(\zeta) = k_M(\zeta) >0$, we have that
$$\overline{\emph{mdim}}_{\emph{M}}(X,G,d) = k(\zeta)\overline{\emph{mdim}}_{\emph{M}}(X,G,\zeta_d)$$
and
$$\underline{\emph{mdim}}_{\emph{M}}(X,G,d) = k(\zeta) \underline{\emph{mdim}}_{\emph{M}}(X,G,\zeta_d).$$
\end{lem}
 \begin{proof}
     From \eqref{mfef}, we have that \begin{align*}
\overline{\text{mdim}}_{\text{M}}(X,G,{d}) & = \limsup_{\varepsilon \to 0} \limsup_{n \to \infty}   \frac{\log r_{F_{n}}(d,\varepsilon,X)}{|F_{n}||\log(\varepsilon)|}\\
& \geq \limsup_{\varepsilon \to 0} \limsup_{n \to \infty}  \frac{\log r_{F_{n}}(\zeta_d,\zeta(\varepsilon),X)}{|F_{n}||\log(\zeta(\varepsilon))|} \frac{|\log(\zeta(\varepsilon))|}{|\log(\varepsilon)|}\\
& = k (\zeta)  \limsup_{\varepsilon \to 0} \limsup_{n \to \infty}  \frac{\log r_{F_{n}}(\zeta_{d},\zeta(\varepsilon),X)}{|F_{n}||\log(\zeta(\varepsilon))|}\\
& = k (\zeta)  \overline{\text{mdim}}_{\text{M}}(X,G,\zeta_d).
\end{align*}
By Proposition \ref{4000} (ii),  we have $$\overline{\text{mdim}}_{\text{M}}(X,G,{d})= k (\zeta)  \overline{\text{mdim}}_{\text{M}}(X,G,\zeta_d).$$

\hspace{4mm}
Similarly, combining  \eqref{mbsgdhd} and Proposition \ref{4000} (i), we can show that \begin{equation*}
\underline{\text{mdim}}_{\text{M}}(X,G,d) = k(\zeta) \underline{\text{mdim}}_{\text{M}}(X,G,\zeta_d). \qedhere
\end{equation*} 
 \end{proof}

\hspace{4mm}
From now on, we will assume that $\rho=\text{diam}_{d}(X)<1.$ Define $${\mathcal{A}}^{+}[0,\rho]:=\{\zeta\in {\mathcal{A}}[0,\rho]: k_m(\zeta)=k_M(\zeta)>0\}.$$  
Note that if $g_{1}, g_{2} \in {\mathcal{A}}^{+}[0, \infty)$, then $g_{1} \circ g_{2} \in {\mathcal{A}}^{+}[0,\infty)$. 
Fix $\zeta \in {\mathcal{A}}^{+}[0,\rho]$, for any $\vartheta \in {\mathcal{A}}^{+}[0,\rho]$ satisfying  $\vartheta(0)=0$, we have   $d(\zeta(x),\vartheta(x))\rightarrow 0$ as $x\rightarrow 0$. For a fixed $\varepsilon >0$, we consider the following set  \begin{equation*}\tilde{B}(\zeta,\varepsilon)=\left\{\vartheta\in {\mathcal{A}}^{+}[0,\rho]:   \zeta(x)(x^{\varepsilon}-1)< \vartheta(x)-\zeta(x) < \zeta(x)\frac{(1-x^{\varepsilon})}{x^{\varepsilon}}, \text{ for  }x\in (0,\rho]\right\}.
\end{equation*}
\hspace{4mm}
Let $\mathcal{T}$ be the topology induced by the sets $\tilde{B}(g,\varepsilon)$, that is, these sets form a subbase for $\mathcal{T}$. The following lemmas come from \cite{ABB}, they are vital to our conclusion.

\begin{lem}\label{4999}\cite{ABB} The map   
\begin{align*}
\mathcal{Z} : ( {\mathcal{A}}^{+}[0,\rho],\mathcal{T})    &\to   (0,1]  \\
g &\mapsto  k(\zeta):=k_{m}(\zeta)
\end{align*}
is continuous.  
\end{lem}

\hspace{4mm}
For the next results, we will consider the set $${\mathcal{A}}^{+}_d(X) =\{\zeta\circ d \in \mathcal{A}_d(X): \zeta\in {\mathcal{A}}^{+}[0,\rho]   \}.$$   

Notice that ${\mathcal{A}}^{+}_d(X)\neq \emptyset$, because  the function  $\zeta(x)=x^a$, for a fixed $a\in (0,1]$, belongs to $\mathcal{A}^{+}[0,\rho]$    (see Example \ref{5001}). In particular, $d\in {\mathcal{A}}^{+}_d(X)$. 

\begin{lem}\cite{ABB}\label{4060}
    Let $X$ be a compact space such that the metric map $d: X\times X \rightarrow [0,\rho]$ is surjective. Then \begin{align*}
\mathcal{Z} : {\mathcal{A}}^{+}[0,\rho]     &\to  {\mathcal{A}}^{+}_{d}(X)   \\
\zeta &\mapsto \zeta \circ d
\end{align*} is a bijective map. 
\end{lem}

\hspace{4mm}
 Suppose that $d: X\times X\rightarrow [0,\rho]$ is surjective. By Lemma \ref{4060}, we can equip  ${\mathcal{A}}^{+}_d(X)$ with  the topology $\mathcal{W}$ such that the map  
\begin{align*}
\mathcal{Z} : ( {\mathcal{A}}^{+}[0,\rho],\mathcal{T})    &\to  ({\mathcal{A}}^{+}_{d}(X),\mathcal{W})  \\
\zeta &\mapsto  d
\end{align*}
 is a homeomorphism. 

\hspace{4mm}
With the above lemmas, we finally have the continuity of upper and lower metric mean dimension in $({\mathcal{A}}^{+}_{d}(X),\mathcal{W}).$
\begin{thm}\label{4061} Let $(X,G)$ be a $G$-system such that the  metric map $d:X\times X \rightarrow [0,\rho]$ is surjective. 
Suppose that $\emph{mdim}_{\emph{M}}(X,G,d)< \infty$. The maps
\begin{align*}
\overline{\emph{mdim}}_{\emph{M}}(X,G,T) \colon ({\mathcal{A}}^{+}_{d}(X),\mathcal{W}) &\to \R  \\
\zeta_{d} &\mapsto \overline{\emph{mdim}}_{\emph{M}}(X,G,\zeta_d)
\end{align*} and \begin{align*}
\underline{\emph{mdim}}_{\emph{M}}(X,G,T) \colon ({\mathcal{A}}^{+}_{d}(X),\mathcal{W}) &\to \R  \\
\zeta_{d} &\mapsto \underline{\emph{mdim}}_{\emph{M}}(X,G,\zeta_{d})
\end{align*}
are continuous.    
\end{thm}
\begin{proof} 
We only prove the case $\overline{\text{mdim}}_{\text{M}}(X,G,T) \colon ({\mathcal{A}}^{+}_{d}(X),\mathcal{W}) \to \R$.

\hspace{4mm}
If $\overline{\text{mdim}}_{\text{M}}(X,G,d)=0$, by Lemma \ref{4004} we know that $\overline{\text{mdim}}_{\text{M}}(X,G,T) \colon {\mathcal{A}}^{+}_{d}(X) \to \R$ is the zero map.  Suppose that $0<\overline{\text{mdim}}_{\text{M}}(X,G,d)<\infty$.  Take $ \tilde{d}$   in  ${\mathcal{A}}^{+}_{d}(X)$ and let $g_{\tilde{d}}$  be the  unique map in ${\mathcal{A}}^{+}[0,\rho]$ such that $ \tilde{d}=\zeta_{\tilde{d}}\circ d $. It follows from Lemma \ref{4004} that $$  \overline{\text{mdim}}_{\text{M}}(X,G,T)(\tilde{d})=\overline{\text{mdim}}_{\text{M}}(X,G,T)(\zeta_{\tilde{d}}\circ d)=\frac{\overline{\text{mdim}}_{\text{M}}(X,G,d)}{k(\zeta_{\tilde{d}})} .$$

 Given that $k(\zeta)>0$ for any $\zeta \in{\mathcal{A}}^{+}[0,\rho]$, then by  Lemma \ref{4999}, $\overline{\text{mdim}}_\text{M}(X,G,T) \colon {\mathcal{A}}^{+}_{d}(X) \to \R$ is continuous.  
\end{proof}

\hspace{4mm}
Finally, we will give some examples of maps $g\in{\mathcal{A}}^{+}[0,\rho]$ and  the respective  expressions for $\text{mdim}_{\text{M}}(X,G,\zeta_d)$. 

\begin{example}\label{5001}   Let $(X,G)$ be a $G$-system with a metric $d$. Fix any $a \in (0, 1]$, we consider the function $\zeta(x) = x^a$ defined for all $x \in [0,\infty)$, it's obvious that $\zeta(x)$ is subadditive. Then define $\zeta_d(x, y) = d(x, y)^a$, it is clear that $k(\zeta) = a$. Therefore, by Lemma \ref{4004}, we have
\begin{equation*}
\underline{\emph{mdim}}_{\emph{M}}(X ,G,\zeta_d) = \frac{\underline{\emph{mdim}}_{\emph{M}}(X,G,d)}{a},~\overline{\emph{mdim}}_{\emph{M}}(X ,G,\zeta_d) = \frac{\overline{\emph{mdim}}_{\emph{M}}(X,G,d)}{a}.
\end{equation*}
By Lemma \ref{3003}, for any F$\phi$lner sequence $\{F_{n}\}$ in $G$, we have that
\begin{equation*}
\underline{\emph{mdim}}_{\emph{H}}(X ,\{F_{n}\},\zeta_{d})=\frac{\underline{\emph{mdim}}_{\emph{H}}(X,\{F_{n}\},d)}{a},~\overline{\emph{mdim}}_{\emph{H}}(X ,\{F_{n}\},\zeta_{d})=\frac{\overline{\emph{mdim}}_{\emph{H}}(X,\{F_{n}\},d)}{a}  .
\end{equation*}
 \end{example}

\begin{example}\label{5003}
Let $(X,G)$ be a $G$-system with a metric $d$. Consider $g(x) = \log(1 + x)$, we have 
that $g(x + y) \leq g(x) + g(y)$.   
Let $g_1(x) = x^a$, for $a \in (0, 1)$, and $g_2(x) = \log (1 + x)$. Then $\vartheta(x)=g_2 \circ g_1(x) = \log (1 + x^a) \in {\mathcal{A}}^{+}[0, \infty)$. Observe that 
$$k (\vartheta) = \liminf_{\varepsilon \to 0^+} \frac{\log(\log(1+\varepsilon^{a}))}{\log(\varepsilon)}=\limsup_{\varepsilon \to 0^+} \frac{\log(\log(1+\varepsilon^{a}))}{\log(\varepsilon)}=a.$$ Therefore, by Lemma \ref{4004}, we know that $$\underline{\emph{mdim}}_{\emph{M}}(X,G,\vartheta_d) = \frac{\underline{\emph{mdim}}_{\emph{M}}(X,G,d)}{a}, ~\overline{\emph{mdim}}_{\emph{M}}(X,G,\vartheta_d) = \frac{\overline{\emph{mdim}}_{\emph{M}}(X,G,d)}{a}.$$

\end{example}

\begin{example} Let $(X,G)$ be a $G$-system with a metric $d$. Suppose that  $\vartheta:X \rightarrow X$ is $ \alpha$-H\"older map for some $\alpha\in(0,1)$,  then there exists $K>0$ such that
\[d(\vartheta(x),\vartheta(y))\leq K d(x,y)^{\alpha}\quad\text{for all }x,y\in X.\]
Define  $ d_\vartheta(x,y) = d(\vartheta(x),\vartheta(y)) $   for all $x,y \in X$, it follows from Example \ref{5001}  that 
$$\underline{\emph{mdim}}_{\emph{M}}(X, G,d_\vartheta) \leq \underline{\emph{mdim}}_{\emph{M}}(X,G, d^{\alpha})=\frac{ \underline{\emph{mdim}}_{\emph{M}}(X,G, d)}{\alpha},$$
$$\overline{\emph{mdim}}_{\emph{M}}(X, G,d_\vartheta) \leq \overline{\emph{mdim}}_{\emph{M}}(X,G, d^{\alpha})=\frac{ \overline{\emph{mdim}}_{\emph{M}}(X,G, d)}{\alpha},$$
and for any F$\phi$lner sequence $\{F_{n}\}$ in $G$, we have that
$$ \underline{\emph{mdim}}_{\emph{H}}(X, \{F_{n}\},d_\vartheta) \leq  \underline{\emph{mdim}}_{\emph{H}}(X,\{F_{n}\}, d^{\alpha})=\frac{ \underline{\emph{mdim}}_{\emph{H}}(X,\{F_{n}\}, d)}{\alpha},$$ 
$$ \overline{\emph{mdim}}_{\emph{H}}(X, \{F_{n}\},d_\vartheta) \leq  \overline{\emph{mdim}}_{\emph{H}}(X,\{F_{n}\}, d^{\alpha})=\frac{ \overline{\emph{mdim}}_{\emph{H}}(X,\{F_{n}\}, d)}{\alpha}.$$  
\end{example} 
 
\section{Acknowledgments}
	 The authors would like to thank their doctoral supervisor, Siming Tu, who provided with much guidance and assistance.

 \scalebox{1.1}{School of Mathematics(Zhuhai), Sun Yat-sen University, }\\
 \scalebox{1.1}{Zhuhai, Guangdong, 519000, P.R. China}
 \end{document}